\documentclass[12pt,twoside]{amsart}

\usepackage{verbatim}
\usepackage{amsmath}
\usepackage{amssymb,latexsym}
\usepackage{enumerate}
\usepackage{amssymb}
\usepackage{a4wide}
\usepackage[latin1]{inputenc}
\usepackage[T1]{fontenc}
\usepackage{times}

\usepackage{amsthm,thmtools,xcolor}

\def\Pr{\mathcal{P}}

\def\hpq0{h^{p,q}_{\leq 0}}
\def\Hpq0{\H_{\leq 0}^{p,q}}

\def\dphi{\partial^{\phi}}
\def\db{\bar\partial}

\def\dbar{\bar\partial}
\def\ddbar{\partial\dbar}

\def\ep{\varepsilon}

\def\C{{\mathbb C}}

\def\U{{\mathcal U}}

\def\cR{{\mathcal R}}

\def\F{{\mathcal F}}
\def\L{{L}}
\def\B{{\mathcal B}}
\def\X{{\mathcal X}}
\def\Y{{\mathcal Y}}
\def\D{\mathcal{D}}
\def\bD{\mathbb{D}}

\def\ol{\overline}
\def\wt{\widetilde}

\def\Proj{{\mathcal Proj}}

\def\K{{\mathcal K}}

\DeclareMathOperator{\Sym}{Sym}

\def\X{{\mathcal X}}
\def\L{{\mathcal L}}
\def\P{{\mathcal P}}
\def\B{{\mathcal B}}
\def\H{{\mathcal H}}
\def\E{{\mathcal E}}

\def\O{{\mathcal O}}

\def\V{{\mathcal V}}

\def\U{{\mathcal U}}

\def\be{\begin{equation}}
\def\ee{\end{equation}}

\newtheorem{thm}{Theorem}[section]
\newtheorem{lma}[thm]{Lemma}
\newtheorem{cor}[thm]{Corollary}
\newtheorem{prop}[thm]{Proposition}

\theoremstyle{definition}

\newtheorem{df}{Definition}

\theoremstyle{remark}

\newtheorem{preremark}{Remark}
\newtheorem{preex}{Example}

\newenvironment{remark}{\begin{preremark}}{\qed\end{preremark}}

\numberwithin{equation}{section}

\begin{document}

\title[]
{Algebraic fiber spaces and \\ curvature of Higher direct images}

\author[]{Bo Berndtsson, Mihai P\u{a}un, Xu Wang}

\address{Bo Berndtsson \endgraf
  Department of Mathematical Sciences,\endgraf
  Chalmers University of Technology and\endgraf
  University of Gothenburg,\endgraf
  SE-412 96 Gothenburg, Sweden\endgraf}
\email{bob@chalmers.se}

\address{Mihai P{\u a}un\endgraf
Mathematisches Institut \endgraf
der Universit\"at Bayreuth\endgraf}
\email{mihai.paun@uni-bayreuth.de}

\address{Xu Wang\endgraf
  DEPARTMENT OF MATHEMATICAL SCIENCES,\endgraf
  NORWEGIAN UNIVERSITY OF SCIENCE AND TECHNOLOGY,\endgraf
  NO-7491 TRONDHEIM, NORWAY\endgraf}
\email{xu.wang@ntnu.no}

\begin{abstract} Let $p:X\to Y$ be an algebraic fiber space, and let $L$ be a line bundle on $X$.
  In this article we obtain a curvature formula
  for the higher direct images of $\Omega^{i}_{X/Y}\otimes L$ restricted to a suitable Zariski open subset of $X$. Our results are particularly
meaningful in case
  $L$ is semi-negatively curved on $X$ and strictly negative or trivial on smooth fibers of $p$. Several applications are obtained, including a new proof of
a result by Viehweg-Zuo in the context of canonically polarized family of maximal variation and its version for Calabi-Yau families. The main feature of our approach is that the general curvature formulas we obtain allow us to bypass the use of ramified covers --and the complications which are induced by them.

\end{abstract}

\bigskip

\maketitle

\section{Introduction}

Let $p: X\to Y$ be a surjective holomorphic map of relative dimension $n$ between two non-singular, compact K\"ahler manifolds. We denote by 
$\Delta\subset Y$ 
the union of codimension one components of the set of singular values of $p$, and we will use the same notation for the associated divisor. The reduced divisor corresponding to the support of the inverse $p^{-1}(\Delta)$ will be denoted by $W$. In order to formulate our
results, we will assume from the very beginning that $\Delta$ and $W$ have simple normal
crossings.
\medskip

\noindent We denote by $\Omega^1_X\langle W\rangle$ and $\Omega^1_Y\langle \Delta\rangle$ the logarithmic co-tangent bundle 
associated to $(X, W)$ and $(Y, \Delta)$ respectively (cf. \cite{Del} and the references therein).
The pull-back of forms with log poles along $\Delta$
gives an injective map of sheaves
\be\label{higher0} p^\star\Omega^1_Y\langle \Delta\rangle\to \Omega^1_X\langle W\rangle.
\ee
Let $\displaystyle \Omega^1_{X/Y}\langle W\rangle$ be the co-kernel
of the map \eqref{higher0}.

\noindent We consider a Hermitian line bundle $(L, h_L)$ on $X$.
The metric $h_L$ is assumed to be smooth when restricted to the generic fibers of $p$. The curvature of $(L, h_L)$ will be denoted by $\Theta_L$.
\smallskip

\noindent 
In this article our primary goal is to explore the differential geometric properties of the higher direct images 
\be\label{higher1}
\cR^{i}p_\star\left(\Omega^{n-i}_{X/Y}\langle W\rangle\otimes L\right)
\ee
where $n\geq i\geq 1$ together with some applications. We will proceed in two steps. 

\medskip 

\noindent (A) \emph{Restriction to the set of regular values}
First we will analyze the properties of \eqref{higher1} when restricted
to the complement of a closed analytic subset of $Y$ which we now describe.
By a fundamental result of Grauert \cite{BaSt}, the sheaves
\eqref{higher1} are coherent, thus 
locally free on a Zariski open subset, $U$,  of $Y$.

The set of singular values
 of the map $p$ will be denoted by $\Sigma\subset Y$.
 Then the complement $\Sigma\setminus\Delta$ is a union of
 algebraic subsets of $Y$ of codimension at least two. For each $i= 1,\dots, n$
we consider the function
\be\label{higher50} t\to h^{i}(X_t, \Omega^{n-i}_{X_t}\otimes L_t),
\ee
defined on $Y_0:= Y\setminus \Sigma$. These functions are
upper semicontinuous, hence constant 
on a Zariski open subset of $Y_0$ (in the analytic topology) which will be denoted by ${\mathcal B}$.

The set $p^{-1}(\B)\subset X$ will be denoted by
$\X$, and we will use the same symbol $p:\X\to \B$ for the restriction to $\X$ of our initial map.
It is a proper holomorphic submersion and the
restriction of \eqref{higher1} to $\B$ is equal to
\be\label{higher2}
\cR^{i}p_\star\big(\Omega^{n-i}_{X/Y}\langle W\rangle\otimes L\big)|_{\B}
=
\cR^{i}p_\star\big(\Omega^{n-i}_{\X/\B}\otimes L\big).
\ee 
 \eqref{higher2} thus defines a vector bundle, say $\mathcal H^{i, n-i}$, over $\B$. Since we have now restricted to the set where \eqref{higher50} is constant, its 
fiber at $t$ is isomorphic to the
Dolbeault cohomology space 
$\displaystyle H^{i}(X_t, \Omega^{n-i}_{X_t}\otimes L_t)$ (see \cite{BaSt}, Theorem 4.12). In section 3 we will give another description of this vector bundle which is more suitable for our computations. 
\smallskip

\noindent Let $\displaystyle [u]\in \Gamma(V, \H^{i, n-i})$ be a local holomorphic section of \eqref{higher2}, where 
$V\subset \B$ is a coordinate open subset and $u$ is a $\dbar$-closed $(0,i)$-form with values in the bundle $\displaystyle \Omega^{n-i}_{\X/\B}\otimes L|_{p^{-1}(V)}$. The restriction 
$u_t$ of $u$ to $X_t$ is a $\dbar$-closed $L_t$-valued $(i,n-i)$-form on $X_t$. We denote by $[u_t]$  its $\dbar$-cohomology class. Let $\omega$ be a smooth $(1,1)$-form on $\X$, whose restriction to the fibers of $p$ is positive definite. Then one may use $\omega|_{X_t}$ and $h_L$ to define a Hermitian norm on $\H^{i,n-i}$. More precisely, the norm of a class $[u_t]$ will be defined as the $L^2$-norm of its unique $\dbar$-harmonic representative, say $u_{t,h}$, i.e.
\be\label{higher3}
\Vert [u]\Vert^2_t:= \int_{\X_t}|u_{t, h}|^2e^{-\varphi}\omega^n/{n!}
\ee
By a fundamental result of Kodaira-Spencer the variation of
$u_{t, h}$ with respect to $t$
is smooth.
\smallskip 

Since the bundle $\Omega^{n-i}_{\X/B}\otimes L$ is a quotient of $\Omega^{n-i}_{\X}\otimes L$, we can construct (cf. section 2) a smooth $(0,i)$-form $\wt u$ with values in
$\Omega^{n-i}_{\X}\otimes L|_{p^{-1}(V)}$ which is $\dbar$-closed on the fibers of $p$ and whose restriction to each $X_t$ belongs to the cohomology class $[u_{t}]$. In what follows $\wt u$
will be
called a \emph{representative} of $[u]$. 
\smallskip

\noindent Summing up, $\H^{i, n-i}$ is a holomorphic Hermitian vector bundle. The sign of the corresponding curvature form is of fundamental importance in applications. Certainly, the curvature tensor is an intrinsic object, completely determined by the complex and the Hermitian structure of the bundle. Nevertheless, it can be expressed
in several ways, depending on the choice of the representatives
$\wt u$ mentioned above. So an important question would be to choose "the best" representative for the holomorphic sections of $\H^{i, n-i}$.
\smallskip


\smallskip

\noindent We first assume that
the curvature form $\Theta_L$ of $(L, h_L)$ is semi-negative
on $X$ and strictly negative on the fibers of $p:\X\to \B$.
In this case, we can take $\omega:= -\Theta_L$. Another consequence of this assumption is the existence of a 
\emph{vertical representative} corresponding to 
 any section $u$ (cf. section 5).
 The vertical representatives are used in the proof of
 our first main result, namely Theorem \ref{curvature}. We show that
the curvature of $\H^{i, n-i}$ evaluated in a direction $u$
can be written as difference of two semi-positive
forms on the base $\B$. 
Roughly speaking, the positive contribution in the expression of the curvature
is given by the contraction (or cup-product) with the Kodaira-Spencer class.
\smallskip

\noindent Rather than stating Theorem \ref{curvature}
here, we mention one of its consequences, which will be important in applications.
Let $\K^i$ be the kernel of the map
\be\label{higher4}
\cR^{i}p_\star\big(\Omega^{n-i}_{X/Y}\langle W\rangle\otimes L\big)
\to \cR^{i+1}p_\star\big(\Omega^{n-i-1}_{X/Y}\langle W\rangle\otimes L\big)\otimes \Omega^1_Y\langle \Delta\rangle
\ee
defined by contraction with the Kodaira-Spencer class. By shrinking $\B$, we can assume that the restriction $\K^i|_\B$ is a sub-bundle of
$\cR^{i}p_\star\big(\Omega^{n-i}_{X/Y}\langle W\rangle\otimes L\big)|_\B$. We have the following statement.

\begin{thm}\label{main, 0} We assume that $\Theta_L\leq 0$ on $X$
and that $\Theta_L|_{\X_t}< 0$ for any $t\in \B$. Then   
the bundle $\displaystyle \K^i|_\B$ endowed with the metric induced from \eqref{higher3} is semi-negative in the sense of Griffiths,
provided that we choose $\omega|_{\X_t}= -\Theta_L|_{\X_t}$.
\end{thm}  
\medskip

We list next a few other results related to Theorem \ref{curvature}. 
For example, assume that the canonical bundle $K_{\X_t}$ of $\X_t$ is ample
(i.e. the family $p$ is canonically polarized).
Then, by a crucial theorem of Aubin and Yau (\cite{Aubin}, \cite{Yau}), we can construct a metric $h_t= e^{-\psi_t}$ on $\displaystyle K_{\X_t}$
such that $\omega_t= \sqrt{-1}\ddbar\psi_t$ is 
K\"ahler-Einstein. By a result of Schumacher, \cite{Schumacher}, the metric on the relative canonical bundle $\displaystyle K_{\X/\B}$ induced by the family
$\left(e^{-\psi_t}\right)_{t\in \B}$ is semi-positively curved. 
We then consider the induced metric on
$\displaystyle L:= -K_{\X/\B}$. Our formula in this case coincides with earlier results of Siu, \cite{Siu}, Schumacher, \cite{Schumacher} and To-Yeung, \cite{To-Yeung}, and following the arguments in \cite{Schumacher} and \cite{To-Yeung} (with some variations), we show that this leads to the construction of a metric of strictly negative sectional curvature on the base manifold $Y$.  

During the (long) preparation of this paper, the article \cite{Naumann} by Naumann was posted on arXiv. In this work the Theorem \ref{curvature} is also established. Naumann shows that the expression of the curvature of $\H^{i, n-i}$ can be derived by a  
method similar to Schumacher's,
\cite{Schumacher}.
Our approach here is very different
and perhaps lighter technically. It can be seen as a generalization of the computations in \cite{Berndtsson1} and \cite{Berndtsson2} for the case of $\H^{n,0}$ (or better say, $\H^{0,n}$).
\medskip

\noindent In section 7 we adapt the previous arguments to the case when the curvature of the bundle $L$ is not strictly negative but instead identically zero on the fibers, and seminegative on the total space. We then establish a variant of Theorem \ref{curvature} in this setting, Theorem 7.5. This formula generalizes the classical formula of Griffiths for the case when $L$ is trivial. As it turns out, the only difference between our formula and Griffiths' is an additional term, coming from the curvature of $L$ on the total space.  In previous  work, Nannicini \cite{Nannicini}, has generalized Siu's computations from \cite{Siu} for $\H^{n-1,1}$ under related conditions when $L$ is the relative canonical bundle. Assuming also that the fibers have trivial canonical bundle, C-L Wang ( \cite{CLWang}) has simplified Nannicini's argument by using Griffiths' theorem instead of Siu's method. Here we follow the same argument as Wang, replacing Griffiths' theorem by Theorem 7.5 to be able to treat also the case when the canonical bundles of the fibers is only assumed to be flat.  We also mention that the general case, when the curvature of $L$ is seminegative on fibers is treated in a somewhat different way by the third author in \cite{Wang-h}.

\medskip

\noindent  A projective family $p:X\to Y$ is called
\emph{Calabi-Yau} if $c_1(X_y)= 0$, which is equivalent to the fact that some multiple of the canonical bundle of the generic fibers is trivial. In this case $L=-K_{\X/\B}$ is fiberwise flat and we show that $L$ has a metric which is moreover seminegative on the total space. We say that $p$ has \emph{maximal variation} if the Kodaira-Spencer map is injective on $Y\setminus \Sigma$.
\smallskip

\noindent The notion of  \emph{Kobayashi hyperbolicity} describes in a very accurate manner weather a complex manifold contains "large" discs. In this setting we have the following statement, obtained as a corollary of our curvature formules combined with a few ideas from
\cite{To-Yeung}, \cite{Schumacher}.

\begin{thm}\label{hypb1}
Let $p:X\to Y$ be a canonically polarized or Calabi-Yau family. We assume that 
$p$ has maximal variation. Then there is no non-constant holomorphic curve $f: \C\to Y
\setminus \Sigma$.
\end{thm}

\noindent 
Our result \ref{hypb1} is a direct consequence of a more general statement
we establish in subsection 8.3.

\noindent We note that Theorem \ref{hypb1} (and a few other results we prove in this article)  have  also been announced by S.-K. Yeung and
W. To in their joint project \cite{Yeung-Ober}.

\bigskip

\noindent (B) \emph{Extension across the singularities.}
 The kernel
\be\label{higher5}
\K^i\subset \cR^{i}p_\star\big(\Omega^{n-i}_{X/Y}\langle W\rangle\otimes L\big)
\ee
of the map \eqref{higher4} is a coherent sheaf of
${\mathcal O}_Y$-modules whose sections are called \emph{quasi-horizontal} in \cite{Griff84} (p 26). Unlike in the case $L= \O_X$,
it is not clear whether $\K^i$ is torsion free or not (cf. \cite{Kang}). Fortunately, this is not a problem for us. We define
\be\label{tquot}
\K^i_f:= \K^i/T(\K^i)
\ee
where $T(\K^i)\subset \K^i$ is the torsion subsheaf. It turns out that the metric induced on $\K^i_f|_\B$ is also semi-negatively curved in the sense of Griffiths (the quotient map $\K^i\to \K^i_f$ is an isometry).
\medskip

\noindent Our goal is to show that the
metric defined on the
bundles $\K^i_f|_\B$ extends as a semi-negatively curved singular Hermitian metric (in the sense of \cite{BP}, \cite{PT}) on the torsion-free sheaf $\K^i_f$, provided that the metric $h_L$ satisfies one of the two conditions below. 
\smallskip

\noindent $\left({\mathcal H}_1\right)$ We have
$\displaystyle \Theta_L\leq 0$ on $X$
and moreover $\displaystyle \Theta_L|_{X_y}= 0$ for each $y$ in the complement of some Zariski closed set.
\smallskip

\noindent $\left({\mathcal H}_2\right)$ We have
$\displaystyle \Theta_L\leq 0$ on $X$
and moreover there exists a K\"ahler metric $\omega_Y$ on $Y$ such that we
have $\displaystyle \Theta_L\wedge p^\star\omega_Y^m\leq -\varepsilon_0\omega\wedge p^\star\omega_Y^m$ on $X$ (where $m$ is the dimension of the base).
\smallskip

\noindent Thus the first condition requires $L$ to be semi-negative and
trivial on fibers, whereas in $\left({\mathcal H}_2\right)$ we assume that $L$
is uniformly strictly negative on fibers, in the sense that  we have 
\be\label{mp0301'} 
\Theta_{h_L}(L)|_{X_y}\leq - \varepsilon_0\omega|_{X_y} 
\ee
for any regular value $y$ of the map $p$.

\medskip

\noindent In this context we have the following result.

\begin{thm}\label{main, I}
  Let $p:X\to Y$ be an algebraic fiber space, and let $(L, h_L)$ be a Hermitian
  line bundle which satisfies one of the hypothesis
  $\left({\mathcal H}_i\right)$ above. We assume that the restriction of $h_L$
  to the generic fiber of $p$ is non-singular. Then: 
  \begin{enumerate}
  
  \item[{\rm (i)}] For each $i\geq 1$ the 
  sheaf $\displaystyle \K^i_f$ admits a semi-negatively curved singular Hermitian metric. 
  \smallskip
  
    \item[{\rm (ii)}] We assume that curvature form of $L$ is smaller than $-\varepsilon_0p^\star\omega_Y$ on the $p$-inverse image of some open subset $\Omega\subset Y$ of $Y$. Then the metric on $\displaystyle \K^i_f$ is strongly negatively curved on $\Omega$ (and semi-negatively curved outside a codimension greater than two subset of $Y$).
\end{enumerate}    
\end{thm}

\noindent As a consequence of the point (i),
we infer e.g. that the dual sheaf $\displaystyle \K^{i\star}_f$ is weakly semi-positive in the sense of Viehweg, cf. \cite{Vbook}, \cite{Kang},\cite{BruneB}.
We note that similar statements   
appear in various contexts in algebraic geometry, cf. \cite{BruneB}, \cite{GT},
\cite{Fujino}, \cite{KK08}, \cite{Kol87}, \cite{MT07},
\cite{MT08}, \cite {PY}, \cite{PoSch}, \cite{VZ02}, \cite{VZ03},
\cite{Kang} among many others.
The fundamental results of Cattani-Kaplan-Schmid \cite{CKS} are an indispensable tool in the arguments of most of the articles quoted above.
Unfortunately, the analogue of the period mapping in our context is not defined (because of the twisting with the bundle
$L$) so we do not have these techniques at our disposal.  
\smallskip

\noindent Part of the motivation for the analysis of the curvature properties of the kernels $\K^i$ is the existence of a non-trivial map
$$\K^{i\star}_f\to \Sym^i\Omega^1_Y\langle \Delta\rangle$$
for some $i$ such that $\displaystyle \K^{i\star}_f\neq 0$ provided that $L^\star:= \det \Omega^1_{X/Y}\langle W\rangle$ and $p$ is either Calabi-Yau or canonically polarized with
maximal variation. This is a well-known and very important fact. We refer to \cite{VZ03}, \cite{PoSch} and the references therein for related results.

\smallskip

\noindent The point is that here we construct the so-called
\emph{Viehweg-Zuo sheaf} $\K^{i\star}$ in a very direct and explicit manner, without the use of ramified covers as in the original approach. 
Indeed, in case of a canonically polarized or Calabi-Yau family $p$ with maximal variation we show in section 9 that the bundle $\det \Omega^1_{X/Y}\langle W\rangle$
has a property which is very similar to the hypothesis $\left({\mathcal H}_i\right)$ above. Hence Theorem \ref{main, I} applies, and 
we obtain a new proof of the existence of the Viehweg-Zuo sheaf, cf. \cite{VZ03}. Our arguments work as well in case of Calabi-Yau families, where we obtain a similar result. We remark that the proof is much simpler than in the canonically polarized case.
\smallskip

On the other hand, in the article \cite{PoSch}, Popa-Schnell obtain 
a vast generalization of the result by Viehweg-Zuo \cite{VZ03} in the case of a family $p$ whose generic fiber is of general type. To this end,
they are using deep results from the theory of Hodge modules.
It would be very interesting to see if our explicit considerations here could provide an alternative argument. 

\medskip

\noindent This paper is organized as follows.

\tableofcontents

\subsection*{Acknowledgments} We would like to thank Christian Schnell,
Valentino Tosatti, St\'ephane Druel and Ya Deng for numerous useful discussions about the topics of this paper. M.P. was partially supported by NSF Grant DMS-1707661 and Marie S. Curie FCFP. It is our privilege to thank the referee for numerous and valuable suggestions.


\section{Preliminaries}

As in the introduction we let  $p:\X\to \B$ be a smooth proper fibration of relative dimension $n$  over an $m$-dimensional base $\B$.  We let $L\to \X$ be a holomorphic line bundle, equipped with a smooth metric, $\phi$. Mostly we will assume that $\phi$ has semi-negative curvature, and that $\Omega:-i\ddbar\phi$ is strictly positive on each fiber $X_t$. Sometimes we let $\omega_t$ denote the restriction of $\Omega$ to the fiber $X_t$, sometimes we just write $\omega$.

We consider maps $t\to u_t$ defined for $t$ in $\B$, where $u_t$ is an $n$ form on $X_t$ for each $t$. We call such a map a {\it section}, of the bundle of $n$-forms on the fiber. We will need to define what it means for such a map to be smooth.

If $t=(t_1,...t_m)$ are local coordinates we can, via the map $p$ consider $t_j$ to be functions on $\X$, and similarily we have the forms $dt_j$, $d\bar t_j$, $dt=dt_1\wedge...dt_m$ and $d\bar t$ on $\X$. That $p$ is a smooth fibration means that its differential is surjective everywhere, so locally near a point $x$ in $\X$, we can complete the $t_j$:s with functions $z_1, ...z_n$ to get local coordinates on $\X$. Then, for $t$ fixed, $z_j$ are local coordinates on $X_t$, and we can write
$$
u_t= \sum_{|I|=p, |J|=q} c_{I,J}(t,z) dz_{I}\wedge d\bar z_J.
$$
We say that $u_t$ is smooth if the coefficients $c_{I,J}$ are smooth. Since this amounts to saying that the form on $\X$
$$
u_t\wedge dt\wedge d\bar t
$$
is smooth, it does not depend on the choice of coordinates.

\begin{lma} A section $u_t$ is smooth if and only if there is a smooth $n$-form, $\tilde u$ on $\X$, such that the restriction of $\tilde u$ to each $X_t$ equals $u_t$.
\end{lma}

\begin{proof} Locally, this is clear from the discussion above: we may take $\tilde u=\sum c_{I,J} dz_I\wedge d\bar z_J$. The global case follows via a partition of unity.  
  \end{proof}
  
These definitions extend to forms with values in $L$. We will call a form $\tilde u$ as in the lemma a {\it representative} of the section $u_t$.
\begin{lma} 
  Two forms, $\tilde u$ and $\tilde u'$,  on $\X$ have the same restriction to the fiber $X_t$, i e represent the same section $u_t$,   if and only if
  $$
  (\tilde u -\tilde u')\wedge dt\wedge d\bar t=0,
  $$
on $X_t$,   where $dt=dt_1\wedge...dt_m$ for some choice of local coordinates on the base. This holds if and only if
  $$
  \tilde u'=\tilde u +\sum_1^m a_j\wedge dt_j+\sum_1^m b_j\wedge d\bar t_j
  $$
  for some $(n-1)$-forms $a_j$ and $b_j$ on $\X$.
\end{lma}
\begin{proof} Both statements follow by choosing local coordinates $(t,z)$ as above.
 \end{proof}

  Notice that the forms $a_j$ and $b_j$ are not uniquely determined, but their restriction to the fibers are uniquely determined.
This follows since if
$$
\sum a_j\wedge dt_j+\sum b_j\wedge d\bar t_j=0,
$$
e g wedging with $\widehat{dt_j}\wedge d\bar t$, we get
$$
a_j\wedge dt\wedge d\bar t=0,
$$
which means that $a_j$ vanishes on fibers by the lemma.

We will mainly be interested in the case when $u_t$ are $\dbar$-closed of a fixed bidegree $(p,q)$ and in their cohomology classes $[u_t]$ in $H^{p,q}(L_t)$. By Hodge theory, each such class has a unique harmonic representative for the $\dbar$-Laplacian $\Box''$ defined by the fiber metric $\phi$ (restricted to $X_t$) and the K\"ahler form $\omega_t$. Since $\Box''=\dbar\dbar^*+\dbar^*\dbar$ and $X_t$ is compact, a form $u$ is harmonic if and only if $\dbar u=0$ and $\dbar^* u=0$. In the next lemma we give an alternative characterization of harmonic forms that will be useful in our computations. 

\begin{lma} Let $X$ be an $n$-dimensional compact complex manifold and $L\to X$ a holomorphic line bundle with a negatively curved smooth  metric $\phi$, i e $\omega:=-i\ddbar\phi>0$ on $X$. Let $u$ be  a $(p,q)$-form on $X$ with values in $L$, where $p+q=n$. Then $u$ is harmonic for the K\"ahler metric $\omega$ on $X$ and the metric $\phi$ on $L$ if and only if $\dbar u=0$ and $\dphi u=0$, where $\dphi=e^\phi\partial e^{-\phi}$ is the $(1,0)$-part of the connection on $L$ induced by the metric $\phi$. Moreover, a harmonic $(p,q)$-form with $p+q=n$ is primitive, i e satisfies $\omega\wedge u=0$.
\end{lma}
\begin{proof} Let $\Box'=\dphi(\dphi)^*+(\dphi)^*\dphi$ be the $(1,0)$-Laplacian. By the Kodaira-Nakano formula
  $$
  \Box''=\Box'+[i\ddbar\phi,\Lambda_\omega].
  $$
  Since $i\ddbar\phi=-\omega$ and the commutator $[\omega,\Lambda_\omega]=0$ on $n$-forms, we see that $\Box'=\Box''$ on $n$-forms. Hence, $\dphi u=0$ if $u$ is harmonic, which proves one direction in the lemma. We also have, when $\dbar u=\dphi u=0$, that $\dbar\dphi u +\dphi\dbar u=0$, which gives that $\omega\wedge u=0$, so $u$ is primitive.

Then $* u=c u$, where $c\in\mathbb C$ such that $|c|=1$, so when $\dphi u=0$, $\dbar^* u=- *\dphi * u=0$. Hence $\Box'' u=0$ if $\dbar u=0$ and $\dphi u=0$, which proves the converse direction of the lemma.
\end{proof}
A similar result of course holds (with the same proof) in case $L$ is positive and we use $\omega=i\ddbar\phi$ as our K\"ahler metric.

At this point we also insert an estimate that will be useful later when we simplify the curvature formula. It generalizes a formula  in \cite{Berndtsson2} (see section 4 in that paper). 
\begin{prop}
  With notation and assumptions as in Lemma 2.3, let $g$ be an $L$-valued form of bidegree $(p,q)$ and $f$ an $L$-valued form of bidegree $(p+1,q-1)$, where $p+q=n$. Assume
  $$
  \dphi g=\dbar f
  $$
  and that $f\perp Range(\dbar)$ and $g\perp Ker(\dphi)$. Then
  $$
  \|f\|^2-\|g\|^2=\langle (\Box''+1)^{-1}f,f\rangle.
  $$
  Moreover,
  $g\in R(\Box''-1)$ and
  $$
   \|f\|^2-\|g\|^2=\langle (\Box''-1)^{-1}g,g\rangle.
  $$
\end{prop}
\begin{proof}
 Since $g$ is orthogonal to the space of harmonic forms, $g$ lies in the domain of $(\Box'')^{-1}$. Since the spectrum of $\Box''$ is discrete, $g$ does also lie in the domain of $(\Box''-\lambda)^{-1}$ for $\lambda$ sufficiently small and positive. We have
  $$
  \lambda\langle (\Box''-\lambda)^{-1} g,g\rangle +\|g\|^2=
  \langle \Box''(\Box''-\lambda)^{-1} g,g\rangle=
  $$
  $$
  \langle \Box'(\Box'-\lambda)^{-1} g,g\rangle =  \langle (\Box'-\lambda)^{-1} \partial^\phi g,\partial^\phi g\rangle=\langle (\Box'-\lambda)^{-1}\dbar f,\dbar f\rangle.
  $$
  Here we have used in the first equality on the second line that $\Box'=\partial^\phi(\partial^\phi)^*+(\partial^\phi)^*\partial^\phi$ and $g$ is orthogonal to $R(\partial^\phi)$. 
  Since $\Box'=\Box''+1$ on forms of degree $(n+1)$, this equals
  $$
  \langle\dbar^*\dbar(\Box''+1-\lambda)^{-1} f,f\rangle=\langle\Box''(\Box''+1-\lambda)^{-1} f,f\rangle,
  $$
  where the last equality follows since $f$ is orthogonal to $R(\dbar)$. 
  Clearly this expression stays bounded as $\lambda$ increases to 1, so the expansion of $g$ in eigenforms  of $\Box''$ can not have any terms corresponding to eigenvalues less than or equal to 1. Hence $g$ lies in the domain of $(\Box''-\lambda)^{-1}$ for all $0\leq\lambda\leq 1$ and
  $$
  \lambda \langle (\Box''-\lambda)^{-1} g,g\rangle +\|g\|^2=
  \langle \Box''(\Box''+1-\lambda)^{-1} f,f\rangle.
  $$
Putting $\lambda=1$ we get
  $$
  \|g\|^2+\langle (\Box''-1)^{-1}g,g\rangle =\|f\|^2,
  $$
  and putting $\lambda=0$ we get
  $$
  \|g\|^2= \langle \Box''(\Box''+1)^{-1} f,f\rangle=\|f\|^2 -
  \langle (\Box''+1)^{-1} f,f\rangle.
  $$
\end{proof}

Note that in the case of  positive curvature, when $i\ddbar\phi=\omega$, we get with the same proof that
$$
\|g\|^2-\|f\|^2=\langle (\Box''+1)^{-1}g,g\rangle,
$$
when $\dbar f=\partial^\phi g$, $g$ is orthogonal to the range of $\partial^\phi$ and $f$ is orthogonal to the kernel of $\dbar$.

For future use in the case of Calabi-Yau families we also record the following counterpart of the previous proposition, which is proved in a similar (but simpler) way.

 \begin{prop} Let $(X,\omega)$ be a compact K\"ahler manifold and $(L,e^{-\phi)}$ a hermitian holomorphic line bundle over $X$ with $i\ddbar\phi=0$.
    Let $g$ be an $L$-valued form of bidegree $(p,q)$ and $f$ $L$-valued of bidegree $(p+1,q-1)$, where $p+q=n$. Assume $f\perp Range(\dbar)$ and $g\perp Ker(\partial^\phi)$ and that
    $$
    \partial^\phi g=\dbar f.
    $$
    Then,
    $$
    \|f\|^2=\|g\|^2.
    $$
  \end{prop}

The next lemma gives in particular a simple way to compute norms of harmonic forms.
\begin{lma} Let $c_n=(i)^{n^2}$. If $u$ is a primitive $(p,q)$-form, with $p+q=n$  on an $n$-dimensional K\"ahler manifold $(X,\omega)$,
  $$
  (-1)^q c_n u\wedge\bar u=|u|^2_\omega \omega^n/n!.
  $$
\end{lma}
The proof can be found in \cite{Griffiths-Harris}.

This means that the $L^2$-norms of primitive forms $u_t$ over a fiber $X_t$  can be expressed in terms of the fiberwise integrals of $u_t\wedge\bar u_t$. When we differentiate such  integrals with respect to $t$ it is convenient to express them as pushforwards.
\begin{lma}
  If $u_t$ is a smooth section of primitive $(p,q)$-forms on $X_t$, where $p+q=n$, we have
  $$
  \int_{X_t}|u_t|^2_{\omega_t} \omega_t^n/n! \,e^{-\phi}=(-1)^q c_n p_*(\tilde u\wedge\overline{\tilde u} e^{-\phi}),
  $$
  if $\tilde u$ is any representative of $u_t$.
  \end{lma}

We will also need (follows from the Fundamental theorem of Kodaira-Spencer in \cite{KSp})

\begin{lma}\label{KodSpe} Let $\X\to \B$ be a smooth proper fibration, and $L\to \X$ a holomorphic line bundle. Assume that the dimension of $H^{p,q}(L_t)$ is a constant. Let $t\to u_t$ be a smooth section of $(p,q)$-forms, $\dbar$-closed on $X_t$, and let $u_{t,h}$ be the harmonic representatives of the classes $[u_t]$ in $H^{p,q}(L_t)$. Then $t\to u_{t, h}$ is smooth.
\end{lma}

In the computations below we will frequently use the {\it interior multiplication} of a form with a vector (field), defined as the adjoint of exterior multiplication under the natural duality between forms and the exterior algebra of tangent vectors. If $V$ is a vector and $\theta$ is a form, it is denoted as  $V\rfloor\theta$. Interior multiplication is an antiderivation, so $V\rfloor (a\wedge b)= (V\rfloor a)\wedge b+ (-1)^{deg(a)}a\wedge (V\rfloor b)$. More generally, we can take the interior product of a wedge product of vectors with a form, and it holds that
  $$
  (V_k\wedge...V_1)\rfloor \theta=V_1\rfloor...(V_k\rfloor \theta).
  $$

\section{The bundles $\H^{n-q,q}$}

We will now consider the vector bundles $\H^{n-q,q}$ over the base $\B$ whose fibers are $$\H^{n-q,q}_t:= H^{n-q,q}(L_t).$$ A smooth section of this bundle can be written as $t\to [u_t]$ such that $t\to u_{t,h}$ (see Lemma \ref{KodSpe}) is a smooth section of the bundle of  $(n-q,q)$-forms on $X_t$ as in Lemma 2.1 (when we write $[u_t]$ it always means that $u_t$ is $\dbar$-closed on $X_t$). We will first give a way to express the $(0,1)$-part of the connection on $\H^{n-q,q}$.

For this we let $\tilde u$ be any representative of the section $[u_t]$, i e an $(n-q,q)$-form on $\X$ whose restriction to each $X_t$ cohomologous to $u_t$. Since $\dbar\tilde u=0$ on fibers, Lemma 2.2 shows that
$$
\dbar\tilde u=\sum_1^m dt_j\wedge  \eta_j +\sum_1^m  d\bar t_j\wedge \nu_j ,
$$
for some forms $\eta_j$ and $\nu_j$. Taking $\dbar$ of this equation and wedging with $\widehat{d\bar t_j}\wedge dt$ we see that $\dbar\nu_j=0$ on fibers.  We can then  define
\be\label{high0}
D''[u_t]=\sum [\nu_j|_{X_t}]\otimes d\bar t_j.
\ee For this to make sense we have to check that it does not depend on the various choices made. First, if $\tilde u'$ is smooth and restricts to each $X_t$ to a form $u_t'$ cohomologous to $u_t$, we can write $u_t'=u_t +\dbar v_t$, where $v_t$ can be taken to depend smoothly on $t$. Then, if $\tilde v$ is a representative of $v_t$ then $\tilde u'-\dbar \tilde v$ is a representative of $u_t$. We therefore only have to check that the definition does not depend on the choice of representatives in the sense of Lemma 2.1. But if we change a representative $\tilde u$ by adding $\sum a_j\wedge dt_j+\sum b_j\wedge d\bar t_j$, $\nu_j$ changes only by the $\dbar$-exact term $\dbar b_j$ so the cohomology class $[\nu_j|_{X_t}]$ does not change. Thus, the definition of $D''[u_t]$ does not depend on the choice of representative of the section $u_t$ nor on the choice of representative in the cohomology class.
\medskip

\begin{prop} The operator $D''$ defines an integrable complex structure on $\H^{n-q,q}$.
\end{prop}
\begin{proof} Let $[u_t]$ be a smooth section and $\tilde u$ a representative of $[u_t]$. Then
  $$
  D''[u_t]=\sum [\nu_j]\otimes d\bar t_j,
  $$
  where
  $$
  \dbar\tilde u=\sum dt_j\wedge \eta_j +\sum d\bar t_j\wedge \nu_j.
  $$
  Since $\dbar\nu_j=0$ on fibers we get
  $$
  \dbar\nu_j=\sum dt_k\wedge\eta_{j k}+\sum  d\bar t_k\wedge \nu_{j k},
  $$
  so
  $$
  D''[\nu_j]=\sum[\nu_{j k}]\otimes d\bar t_k.
  $$
  Using that $\dbar^2\tilde u=0$ we get that $ \nu_{j k}=\nu_{k j}$ on fibers, which gives $(D'')^2[u_t]=0$.
\end{proof}
We now return to the setting in the introduction and let $p:X\to Y$ be a  surjective holomorphic map.
We  take   $\B$ to be the Zariski open subset of $Y$, defined in the introduction (immediately after (1.3)). Then $p$ is a proper submersion over $\B$ and  the $q$:th direct image
\be\label{high1}
\cR^qp_\star(\Omega^{n-q}_{X/Y}\otimes L)|_{\B}
\ee
is the sheaf of sections of a vector bundle with fibers
$\displaystyle H^{q}(X_t, \Omega^{n-q}_{X_t}\otimes L_t)$. These spaces are naturally identified with the spaces $H^{n-q,q}(X_t, L_t)$; the fibers of the bundle $\H^{n-q,q}$.  
We will next  show that 
the holomorphic structure introduced above on $\H^{n-q,q}$ coincides with the usual complex structure of the higher direct image.

By definition, the \emph{holomorphic} sections of the vector bundle \eqref{high1} on a coordinate open subset $V\subset \B$ are elements of
\be\label{high2}
H^{q}\left(p^{-1}(V), \Omega^{n-q}_{\X/\B}\otimes L|_{p^{-1}(V)}\right).
\ee
\smallskip

\noindent 
Let $(u_j)$ 
be a set of $\dbar$-closed $(0, q)$--forms with values in $\Omega^{n-q}_{\X/\B}\otimes L|_{p^{-1}(V)}$
such that their corresponding Dolbeault cohomology classes $[u_j]$ gives a local frame for \eqref{high1} over $V$. Then for each 
$t\in V$ the restriction $\displaystyle(u_j|_{X_t})$
induces a basis for $\displaystyle H^{q}(X_t, \Omega^{n-q}_{X_t}\otimes L)$. In particular we infer that the space of smooth sections of \eqref{high1} coincides with the space of smooth sections of $\H^{n-q,q}$ (as defined at the beginning of the current section), since we can express $[u_t]$ as linear combination of 
$\displaystyle(u_j|_{X_t})$ for $t\in V$.

\begin{lma} 
Under the above identifications, the holomorphic structure on the bundle $\H^{n-q,q}$ coincides with the holomorphic structure on the $q^{\rm th}$ direct image bundle.
\end{lma}
\begin{proof}
  Let $[u]$ be one of the elements in the frame $([u_j])$. By restriction to fibers it defines a smooth section $[u_t]$ of $\H^{q,n-q}$.
We first show that 
\be\label{high3}
D^{\prime\prime}([u_t])=0.
\ee
This can be seen as follows. By a partition of unity we can construct a smooth $(n-q, q)$-form $\wt u$ with values in $L$ whose image by the quotient map 
\be\label{high4}
\Omega^{n-q}_{\X}\otimes L|_{p^{-1}(V)}\to \Omega^{n-q}_{\X/\B}\otimes L|_{p^{-1}(V)}
\ee
is equal to $u$. This implies that $\wt u$ is a representative of $[u]$, and
moreover we have
\be\label{high5}
\wt u\wedge dt= u\wedge dt.
\ee
Notice that in \eqref{high5} the expression $u\wedge dt$
is a well-defined $(n-q+m, q)$-form with values in $L$. Now we apply  $\dbar$ to \eqref{high5} and get
\be\label{high6}
\dbar \wt u\wedge dt= 0
\ee
since $\dbar u\wedge dt$ equals zero. Thus the $\nu_j$ defined by the representative $\tilde u$ are even zero, so $D''[u_t]=0$. Hence any holomorphic section of the $q^{\rm th}$ direct image bundle is holomorphic as a section of $\H^{n-q,q}$. 
\smallskip
  
\noindent Let now $[u_t]$ be a smooth section of $\H^{n-q,q}$ such that
\eqref{high3} holds at each point of $V$. We use the  holomorphic frame
$[u_j]$ of $\cR^qp_\star(\Omega^{n-q}_{X/Y}\otimes L)|_{\B}$. We have already shown that the $[(u_j)_t]$ satisfy $D''[(u_j)_t]=0$, so they form a holomorphic frame for the bundle $\H^{n-q,q}$ as well.   Then 
\be\label{high20}
[u_t]= \sum_j f_j[(u_j)_t]
\ee
for some smooth  functions $f_j$ defined on $V$. We have
\be\label{high21}
0=D^{\prime\prime}([u_t])= \sum_j [(u_j)_t]\otimes\dbar f_j,
\ee
so $f_j$ are holomorphic. Then clearly $[u_t]$ is holomorphic as a section of $
\cR^qp_\star(\Omega^{n-q}_{X/Y}\otimes L)|_{\B}$ as well.
\end{proof}

\medskip

\noindent We next equip $\H^{p,q}$ with the hermitian metric
$$
\|[u_t]\|_t^2=\int_{X_t} |u_{t,h}|^2_\omega e^{-\phi} \omega^n/n!
$$
(where we have written $\omega$ instead of $\omega_t$), where $u_{t,h}$ is the harmonic representative of $[u_t]$. 
By Lemmas 2.3 and 2.5, this equals
$$
(-1)^q c_n\int_{X_t} u_{t,h}\wedge\bar u_{t,h}e^{-\phi},
$$
which by Lemma 2.6 also can be written as
\be
\|[u_t]\|_t^2=(-1)^q c_n p_*(\tilde u\wedge\overline{\tilde u} e^{-\phi}),
\ee
for any choice of representative of the smooth section $u_{t,h}$. 

Let now $\tilde u$ be any representative of a smooth section $u_t$ such that each $u_t$ is harmonic on $X_t$. By Lemma 2.3, this means that $\dbar u_t=0$ on $X_t$ and $\dphi u_t=0$ on $X_t$. We then have that
$$
\dphi \tilde u=\sum dt_j  \wedge \mu_j +\sum d\bar t_j\wedge  \xi_j,
$$
by Lemma 2.2, since the left hand side vanishes on fibers. We thus have two sets of equations for $\tilde u$,
\be
\dbar\tilde u=\sum_1^m dt_j\wedge  \eta_j +\sum_1^m  d\bar t_j\wedge \nu_j ,
\ee
and
\be
\dphi \tilde u=\sum dt_j  \wedge \mu_j +\sum d\bar t_j\wedge  \xi_j.
\ee
Taking $\dbar$ of the first equation and wedging with appropriate forms on the base, we see that $\eta_j$ and $\nu_j$ are $\dbar$-closed on fibers, and taking $\dphi$ of the second we see that 
$\mu_j$ and $\xi_j$ are $\dphi$-closed. We stress that the forms $\eta_j$ etc, and even their restrictions to fibers, depend on the choice of representative $\tilde u$, but the cohomology classes they define do not. 

We now polarize (3.10) and use it to compute
$$
\dbar\langle u_t,v_t\rangle_t=(-1)^qc_n \left(\sum p_*(d\bar t_j\wedge\nu_j \wedge\overline{\tilde v} e^{-\phi}) +\sum (-1)^n p_*(\tilde u\wedge\overline{dt_j\wedge  \mu_j} e^{-\phi}) \right )
$$
(the $\eta_j$:s and the $\xi_j$:s give no contrbution for bidegree reasons).

From this we see that $\sum [\mu_{j}]\otimes dt_j=D'[u_t]$, so we now also have an expression for  the $(1,0)$-part of the connection on $\mathcal H^{p,q}$ in terms of representatives.

To compute the curvature of the connection $D=D'+D''$ we will compute
$i\ddbar\|u_t\|_t^2$ and use the classical formula
$$
-\langle \Theta^E_{W,\bar W} u_t, u_t\rangle_t +\|D'_W u_t\|^2_t=\bar W \rfloor  (W \rfloor \ddbar \|u_t\|^2_t),
$$
valid for any holomorphic section of a hermitian vector bundle $E$ and any $(1,0)$ vector in the base, where
\begin{equation}\label{eq:theta-E}
\Theta^E_{W,\bar W}:= \bar W \rfloor  (W \rfloor \Theta^E), \ D'_W:= W\rfloor D',
\end{equation}
and the Chern curvature $\Theta^E$ denotes the square of the Chern connection of the hermitian vector bundle $E$ (in our case $E=\mathcal H^{p,q}$). Here it will be convenient to choose our representatives in a special way, so we first discuss this issue in the next sections.

\section{The horizontal lift.}

Recall that we assume that $\Omega=-i\ddbar\phi$ is strictly positive on fibers. In this situation we can, following Schumacher, \cite{Schumacher}, define the {\it horizontal lift}  of a vector field on the base. This represents a variation of the notion of {\it harmonic lift}, introduced by Siu in \cite{Siu}.
\begin{prop} Let $W$ be a $(1,0)$ vector field on the base $B$. If  $\Omega>0$ on fibers, there is a unique vector field $V=V_W$ on $\X$ such that
  \be
  dp(V)=W
  \ee
  and
  \be
  \Omega(V,\bar U)=0,
  \ee
  for any  field $U$ on $\X$ that is vertical, i e satisfies $dp(U)=0$.
  $V$ is called the horizontal lift of $W$.
  \end{prop}
\begin{proof}
  We begin with uniqueness. Assume $V$ and $V'$ both satisfy (4.1) and (4.2). Then $V-V'$ is vertical, so
  $$
  \Omega(V-V',\overline{V-V'})=\Omega(V,\overline{V-V'})- \Omega(V',\overline{V-V'})=0.
  $$
  Since we have assumed that $\Omega$ is strictly positive on fibers, $V-V'=0$.

  For the existence we note that the proposition is a purely pointwise, linear algebra statement. Given an arbitrary lift $V'$ at a point $x$ in $\X$, a general lift can be written $V=V'+U$, where $U$ is vertical. We therefore have $n$ equations for $n$ unknowns, so uniqueness implies existence.
\end{proof}

Expressed differently, the condition (4.2) means that
$$
V\rfloor \Omega(\bar U)=0
$$
for all vertical $U$, so

\be
V_W\rfloor\Omega =\sum b_j d\bar t_j=:b,
\ee
for some functions $b_j$ if $t_j$ are local coordinates on the base. Clearly $b$ depends linearlily on $W$, so if $W=\sum W_j\partial/\partial t_j$,
$$
V_W\rfloor\Omega=i\sum c_{j,k}(\Omega) W_j d\bar t_k= V_W\rfloor C(\Omega),
$$
where
$$
C(\Omega)=i\sum c_{j,k}(\Omega) dt_j\wedge d\bar t_k.
$$
(Here we have also used that $V_W\rfloor dt_j=W_j$.) We shall see below that $C(\Omega)$ is globally well defined.
We say that a $(1,0)$ vector (field) on $\X$ is horizontal if it is the horizontal lift of some vector (field) on the base. The dimension of the space of horizontal vectors at a point is the dimension of the base, and the intersection between the horizontal vectors and the vertical vectors is just the zero vector. Hence any vector at a point in $\X$ has a unique decomposition as a sum of a vertical and a horizontal vector. Similarily, a $(0,1)$ vector (field) $V$ is horizontal respectively vertical if $\bar V$ is horizontal or vertical.  We next introduce the corresponding notions for forms.

\begin{df} A form $\theta$ is {\it horizontal} if $U\rfloor \theta=0$ for any vertical $U$, and  $\theta$ is {\it vertical} if $V\rfloor \theta=0$ for any horizontal vector (field) $V$.
\end{df} 

To make this a little bit more concrete, let $V_j$ be the horizontal lift of $\partial/\partial t_j$, where $t_j$ are local coordinates on the base for $j=1, ...m$. Let $W_k$ for $k=1,...n$ be a system of vertical $(1,0)$-fields such that $V_j$ and $W_k$ together form a basis for the space of $(1,0)$ vectors near a given point $x$ in $\X$. Then we may find local coordinates $(s,z)$ near $x$ such that $W_k=\partial/\partial z_k$ for $k=1,...n$ and $V_j=\partial/\partial s_j$ for $j=1,...m$ at $x$. (Then, still at $x$, $ds_j(V_k)=\delta_{j k}$ and $ds_j(W_k)=0$, from which it follows that $ds_j=p^*(dt_j)(=dt_j)$.) Then $V_j\rfloor dz_k=0$ so all $dz_k$ are vertical forms at $x$, and similarily $W_k\rfloor ds_j=0$, so $dt_j=ds_j$ are horizontal (which is easy to see  directly). Therefore the horizontal forms are of the type
$$
\sum a_{J,K} dt_J\wedge d\bar t_K,
$$
and vertical forms are of the type
$$
\sum b_{J,K} dz_J\wedge d\bar z_K.
$$

In this terminology, $C(\Omega)$ is horizontal, and the relation $V_W\rfloor\Omega=V_W\rfloor C(\Omega)$ for all horizontal $W$ implies that $\Omega=C(\Omega) + \Omega_v$, where $\Omega_v$  is vertical. This implies that $C(\Omega)$ does not depend on the choice of local coordinates. Indeed, if $C'(\Omega)$ is defined using different local coordinates, $C(\Omega)-C'(\Omega)$ is both horizontal and vertical, hence zero. 

Next, let $V_1, V_2, ...$ be the horizontal lifts of $\partial/\partial t_1, \partial/\partial t_2, ..$, where $t_j$ are local coordinates on the base. Let $\V=V_m\wedge... V_1$ and put for an arbitrary form $\theta$ on $\X$,
\be
P_v(\theta)=(\bar\V\wedge \V)\rfloor(dt\wedge d\bar t\wedge \theta).
\ee
Clearly, $P_v(\theta)$ is always vertical and does not depend on the choice of local coordinates. Moreover, if $\theta$ is already vertical, then $P_v(\theta)=\theta$, so $P_v$ is a projection from the space of all forms  on $\X$ to the space of vertical forms.
\begin{prop} If $\theta$ is an arbitrary form on $\X$, $P_v(\theta)-\theta$ vanishes on fibers. If a vertical form vanishes on fibers it vanishes identically.
\end{prop}
\begin{proof} The first statement means that
  $$
  P_v(\theta)\wedge dt\wedge d\bar t=\theta\wedge dt\wedge d\bar t.
  $$
  This follows from the definition of $P_v$ if we use that contraction is an antiderivation and $ (\bar\V\wedge\V) \rfloor (dt\wedge d\bar t)=1$.

  For the second statement we use that if $\theta$ is vertical, $P_v(\theta)=\theta$. On the other hand, $P_v(\theta)=0$ if $\theta$ vanishes on fibers, since $dt\wedge d\bar t\wedge \theta=0$ then.
  \end{proof}

Finally we note that since $C(\Omega)^{m+1}=0$ and $\Omega_v^{n+1}=0$,
$$
\Omega^{n+m}/(n+m)!=C(\Omega)^m/m!\wedge \Omega_v^n/n!=C(\Omega)^m/m!\wedge\Omega^n/n!.
$$
Thus, if $m=1$, and we write $C(\Omega)=c(\Omega)idt\wedge d\bar t$ we have:

\begin{prop} If $m=1$, 
  $$
  c(\Omega)=\frac{\Omega^{n+1}/(n+1)!}{\Omega^n/n!\wedge idt\wedge d\bar t}.
  $$
\end{prop}
This is the well known geodesic curvature of the curve of metrics $\phi_t$ in Mabuchi space (if $\phi_t$ depends only on the real part of $t$, and by extension in general).
\section{Vertical representatives}

Let $[u]$ be a section of the bundle of smooth $n$-forms on the fibers. Recall that this means that for each $t$ in the base, $u_t$ is an  $n$-form on the fiber $X_t$, and there is a smooth $n$-form $\tilde u$ on the total space $\X$ (a representative of $u_t$)  which restricts to $u_t$ on each $X_t$.
 \begin{prop} Any smooth section $u_t$ has a unique  vertical representative. 
 \end{prop}
 \begin{proof}
   Let $\tilde u$ be an arbitrary representative and let $\hat u:=P_v(\tilde u)$. Then Proposition 4.2 implies that $\hat u$ is also a representative of $u_t$ (since $\hat u-\tilde u$ vanishes on fibers), and $\hat u$ is by definition vertical. Uniqueness also follows from Proposition 4.2.
   \end{proof}

In the next proposition we shall consider sections that are primitive on fibers, i. e. are such that on each $X_t$, $u_t\wedge \Omega =0$. In terms of representatives, this means that
$$
\Omega\wedge\tilde u\wedge dt\wedge d\bar t=0
$$
if $dt=dt_1\wedge .. dt_m$ for a system  of local coordinates on the base.

  \begin{prop} If $\hat u$ is a vertical representative of a section that is primitive on fibers, then
    $$
    \Omega\wedge\hat u=C(\Omega)\wedge\hat u.
    $$
  \end{prop}
  \begin{proof}
    The form
    $$
    \theta:=(\Omega-C(\Omega))\wedge\hat u
    $$
    is a product of vertical forms, hence vertical. Since $C(\Omega)$ vanishes on fibers and $\hat u$ is primitive on fibers, $\theta$ vanishes on fibers. Hence $\theta=0$ by Proposition 4.2.
  \end{proof}
  
  Notice that this means in particular that $\Omega\wedge\hat u\wedge dt=\Omega\wedge\hat u\wedge d\bar t=0$ since $C(\Omega)\wedge dt=C(\Omega)\wedge d\bar t=0$.

  \begin{prop}
    Let $\hat u$ be a vertical representative of a section $u_t$ that satisfies the equation
    $$
    (\dbar +\partial^\phi)u_t=0
    $$
    on each fiber $X_t$. Let $\D:=(\dbar+\partial^\phi)$ on the total space $\X$. Then
    $$
    \D\hat u=\sum a_j\wedge dt_j+ \sum b_j\wedge d\bar t_j,
    $$
    where $a_j$ and $b_j$ are primitive and satisfy $(\dbar+\partial^\phi)a_j=(\dbar+\partial^\phi) b_j=0$ on fibers.
  \end{prop}
  \begin{proof} Since $(\dbar +\partial^\phi) \hat u$ vanishes on fibers, 
    $$
     \D\hat u=\sum a_j\wedge dt_j+ \sum b_j\wedge d\bar t_j,
     $$
     for some $n$-forms $a_j$ and $b_j$. Since $\D^2=\Omega\wedge$ we get
     $$
     \Omega\wedge\hat u=\sum \D a_j\wedge dt_j+\sum \D b_j\wedge d\bar t_j.
     $$
     Wedging with $\widehat{dt_j}\wedge d\bar t$ we find that  the restriction of $\D a_j$ to fibers vanishes since $\hat u\wedge\Omega\wedge d\bar t=0$, and similarily we find that $\D b_j$ vanishes on fibers. Hence $\D^2 a_j=\D^2 b_j$ also vanish when restricted to fibers, which gives that $a_j$ and $b_j$ are primitive. 
  \end{proof}

  Notice the particular case when $u_t$ has pure bidegree $(p,q)$. Then the assumption means that $\dbar \hat u=\partial^\phi \hat u=0$ on fibers, i. e. that $u_t$ is harmonic on each fiber. The first conclusion  is then that
\be
\dbar\hat u=\sum_1^m dt_j\wedge  \eta_j +\sum_1^m  d\bar t_j\wedge \nu_j ,
\ee
and
\be
\dphi \hat u=\sum dt_j  \wedge \mu_j +\sum d\bar t_j\wedge  \xi_j.
\ee
  where all forms $\eta_j$, $\nu_j$, $\mu_j$ and $\xi_j$ are primitive. This follows from Proposition 5.3 since e g $(-1)^n a_j=\eta_j+\mu_j$ and $\eta_j$ and $\mu_j$  have different bidegrees so each of them must be primitive.

We also get
  \be
  \partial^\phi\eta_j=-\dbar\mu_j
  \ee
  and
  \be
  \partial^\phi\nu_j=-\dbar\xi_j
  \ee
  from $(\dbar+\partial^\phi) a_j=0$ and $(\dbar+\partial^\phi) b_j=0$.

\subsection{Vertical representatives and the Kodaira-Spencer tensor.}
  
We first recall the definition of the Kodaira-Spencer class.
 Let $W$ be a holomorphic $(1,0)$-vector field on $\B$ and let $V$ be a smooth  lift of $W$ to $\X$ (so that $dp(V)=W$). In general, unless the fibration is locally trivial, we can not find a holomorphic lift, and the Kodaira-Spencer class is an obstruction to this. First we define
$$
\kappa_V=\dbar V,
$$
if $V$ is any smooth lift of $W$. Then $\kappa_V$ is a $(0,1)$-form on $\X$ with coefficients in $T^{(1,0)}(\X)$.  Since $dp(V)=W$ is holomorphic, it actually takes values in the subbundle, $F=T(\X/\B)$,  of $T^{(1,0)}$ of vertical vectors, i e vectors $V'$ such that $dp(V')=0$. Any other lift, say $V'$, of $W$ can be written  $V+V''$ where $V''$ is vertical and
$$
\kappa_{V'}=\kappa_V+\dbar V'',
$$
so the cohomology class  of $\kappa_V$ in $H^{(0,1)}(\X, F)$ is well defined. We will call it $[\kappa_W]$. Similarily, we let $[\kappa_W^t]$ be the cohomology class of $\kappa_V$ restricted to $X_t$.

In the sequel we will assume  that the lift of $W$ is taken as the horizontal lift. If $W=\partial/\partial t_j$ for a given system of coordinates on the base we will sometimes write  $V_j$ for the horizontal lift of $\partial/\partial t_j$ and just $\kappa_j$ for $\kappa_{V_j}$.  This depends of course on $\Omega$, so it would be more proper to write  $\kappa_j(\Omega)$, but we will use the lighter notation instead.

Locally, using local coordinates $x=(t,z)$ such that $p(x)=t$, $\kappa_V$ can be written
$$
\kappa_V=\sum_1^n Z_j\otimes d\bar z_j +\sum_1^m T_k \otimes d\bar t_k,
$$
where $Z_j$ and $T_k$ are vector fields tangential to fibers.  Its restriction to a fiber $X_t$ is
$$
\kappa^t_V=\sum_1^n Z_j\otimes d\bar z_j,
$$
where we interpret $z_j$ as local coordinates on $X_t$. If  $\theta$ is a form on $\X$, $\kappa_V$ operates on $\theta$ by letting the vector part of $\kappa_V$ operate by contraction, followed by  taking the wedge product with the form part. The result is called $\kappa_V\cup \theta$, which we sometimes write $\kappa_V.\theta$. In the same way, $\kappa^t_V$ operates on forms on a fiber $X_t$. Notice that since the vector parts of the Kodaira-Spencer forms are vertical, the cup product commutes with restriction to fibers: $(\kappa_V\cup \theta)|_{X_t}= \kappa^t_V\cup (\theta|_{X_t})$.

\begin{prop}
  Let $u_t$ be a smooth section of $(p,q)$-forms such that $u_t$ is harmonic on each fiber $X_t$. Let $\eta_j$ and $\xi_j$ be defined as in (5.1) and (5.2), where $\tilde u=\hat u$ is the vertical representative. Then, on each fiber
  $$
  \eta_j=\kappa^t_j\cup u_t
  $$
  and
  $$
  \xi_j=\overline{\kappa^t_j}\cup u_t.
  $$
 
\end{prop}

\begin{proof}
 In the proofs here we will use the easily verified formulas $\dbar(V\rfloor\theta)=(\dbar V).\theta -V\rfloor(\dbar\theta)$ and $\partial^\phi(\bar V\rfloor\theta)=(\partial\bar V).\theta-\bar V\rfloor (\partial^\phi\theta)$ for any form $\theta$ if $V$ is of type $(1,0)$. (Notice however that there are no similar formulas when $V$ is of type $(0,1)$.)

 We  have $V_j\rfloor \hat u=0$, so taking $\dbar$ we get
  $$
  (\dbar V_j).\hat u=V_j\rfloor\dbar\hat u=\eta_j +...
  $$
  where the dots indicate forms that contain $dt_k$ or $d\bar t_j$ that therefore vanish when restricted to fibers. Restricting to fibers we then get (in view of the remarks above) that
  $$
  \kappa^t_j\cup u_t=\eta_j.
  $$
  For the second statement we use that $\bar V_j\rfloor \hat u=0$ and take $\partial^\phi$. Then
  $$
  (\partial\bar V_j)\rfloor\hat u=\bar V_j\rfloor\partial^\phi\hat u= \xi_j +...,
  $$
  which gives that $\xi_j=\overline{\kappa_j^t}\cup\hat u$.

  \end{proof}

The final result in this section is a reflection of  the familiar fact that the Kodaira-Spencer class defines a holomorphic section of the bundle with fibers
$H^{0,1}(T^{1,0}(X_t))$.

\begin{prop}
  Let $[u_t]$ be a holomorphic section of the bundle $\H^{n-q,q}$. Let $W$ be a holomorphic vector field on $B$. Then
  $$
  [\kappa^t_W\cup u_t]
  $$
  is a holomorphic section of $\H^{n-q-1,q+1}$.
\end{prop}

\begin{proof} It is enough to prove this when $W=\partial/\partial t_j$, so that $\kappa^t_j\cup u_t=\eta_j$ by Proposition 5.4. That $u_t$ is holomorphic means that we can find some representative $\tilde u$ such that
  $$
  \dbar\tilde u=\sum dt_j\wedge\eta_j +\sum d\bar t_j \wedge \nu_j,
  $$
and each $\nu_j$ is zero on fibres, i.e. we can write
$$
\nu_j=\sum dt_k \wedge a_{j}^k +\sum d\bar t_k \wedge b_j^k.
$$
Taking $\dbar$ of $\dbar \tilde u$, we get that
$$
0= \sum dt_j\wedge \dbar \eta_j +\sum d\bar t_j \wedge \dbar \nu_j.
$$
Thus each $\dbar \eta_j$ is zero on fibres and we can write
$$
\dbar \eta_j=\sum dt_k \wedge \eta_j^k +\sum d\bar t_k \wedge \nu_j^k. 
$$
Thus we have
\begin{eqnarray*}
0 & = & \sum dt_j \wedge dt_k \wedge \eta_j^k + \sum  dt_j \wedge d\bar t_k \wedge \nu_j^k \\
& & +\sum dt_k \wedge d\bar t_j \wedge \dbar  a_{j}^k + \sum d\bar t_k \wedge d\bar t_j \wedge \dbar b_j^k,
\end{eqnarray*}
which implies that $\nu_j^k=-\dbar a_k^j$ on fibres. Thus each $\nu_j^k$ are $\dbar$-exact on fibres and thus $\eta_j$ defines a holomorphic section of $\H^{n-q-1,q+1}$.
\end{proof}

As an example of this we consider the case when $L=-K_{\X/B}$. Then $\H^{n,0}$ has a canonical trivializing section, $u^0_t$, defined as follows. On any complex manifold $X$ the bundle $K_X-K_X$ is of course trivial and  has a canonical trivializing section defined locally as $s=dz\otimes (dz)^{-1}$. Applying this to each fiber $X_t$ we get a trivializing section $u^0_t$ of the line bundle $\H^{n,0}$. Wedging with the corresponding trivializing section of the base, which is locally $dt\otimes (dt)^{-1}$ we get an $(n+m,0)$-form, $U^0$ on $\X$ with values in $K_{\X/\B}^{-1}\otimes K_B^{-1}=K_X^{-1}$ (by the adjunction isomorphism). The form $U^0$ is easily seen to be the canonical trivializing section of $K_\X-K_\X$; in particular it is holomorphic. This means that $\dbar u^0\wedge dt=0$, so $\dbar u^0=\sum dt_j\wedge\eta_j$, so $u^0$ is a holomorphic section of $\H^{n,0}$.

Applying the proposition iteratively we get holomorphic sections of $\H^{n-q,q}$ as
$$
\kappa^t_{i_1}\cup...\kappa^t_{i_q}\cup u^0_t
$$
for any multi-index $I=(i_1, ...i_q)$.

\section{The curvature formulas}

We now have all the ingredients to begin computing the  curvature of $\H^{n-q,q}$ with the $L^2$-metric. Let $u_t$ be a smooth section of this bundle, which now means that $\dbar u_t=0$ on $X_t$ for each $t$ and $\partial^\phi u_t=0$ on $X_t$ for each $t$. We take $\hat u$ to be a vertical representative of the section $u_t$ so that
\be
\dbar\hat u=\sum dt_j\wedge\eta_j +\sum d\bar t_j\wedge \nu_j
\ee
and
\be
\partial^\phi\hat u=\sum dt_j\wedge\mu_j +\sum d\bar t_j\wedge\xi_j.
\ee
Then $\nu_j$ and $\mu_j$ have bidegree $(p,q)$, $\eta_j$ have bidegree $(p-1,q+1)$ and $\xi_j$ have bidegree $(p+1,q-1)$ (forms of negative degree should be interpreted as zero). We decompose each of these forms on each fiber into one harmonic part and one part which is orthogonal to harmonic forms, so that e.g.
$(\eta_j)_h$ is the harmonic part of $\eta_j$ and $(\mu_j)_\perp$ is the part of $\mu_j$ which is orthogonal to harmonic forms.

Our main  curvature formula is as follows.

\begin{thm}\label{curvature} 

 Let $L\to \X$ be a line bundle with metric $e^{-\phi}$ where $i\ddbar\phi<0$ on fibers and define the hermitian structure on $\H^{n-q,q}$ using the K\"ahler forms $\omega_t=-i\ddbar\phi|_{X_t}$ on fibers. Then
$$
 \langle\Theta_{\partial/\partial t_j,\partial/\partial \bar t_k}\, u_t,u_t\rangle_t  =
 $$
 $$
 -\langle (\Box''+1)^{-1}(\mu_j)_\perp, (\mu_k)_\perp\rangle_t -\langle (\Box'' +1)^{-1}(\xi_k),(\xi_j)\rangle_t 
  +\langle(\eta_j)_h,(\eta_k)_h\rangle_t -c_n(-1)^q\int_{X_t} c_{j,k}(\Omega)u_t\wedge \bar u_t e^{-\phi},
  $$
  
where $\mu,\xi,\eta$ are defined  in formulas (5.1) and (5.2), using a vertical representative of $u_t$. 
\end{thm}

Recall that by Proposition 5.4, $\eta_j=\kappa_j\cup u_t$, and that $(\eta_j)_h$ is the harmonic representative of $\eta_j$. Similarily, $\xi_j=\overline{\kappa_j}\cup u_t$. As for the term containing $\mu_j$, we can use the second part of Proposition 2.4 and write
$$
\langle (\Box''+1)^{-1} (\mu_j)_\perp,(\mu_k)_\perp\rangle= \langle (\Box''-1)^{-1}(\eta_j)_\perp,(\eta_k)_\perp\rangle,
$$
in this way we see that the whole formula for the curvature can be expressed in terms of $u_t$ and $\kappa$. 

Notice that all the terms in the curvature formula are negative, except the one coming from the harmonic part of $\eta$. This 'bad term' is clearly zero when $(p,q)=(0,n)$, for bidegree reasons. It also vanishes  in case the harmonic part of $\kappa$ vanishes, so that the fibration is isotrivial. More generally, if $W=\sum a_j\partial/\partial t_j$ is a vector at a point $t$ in the base, and the cohomology class defined by $[\kappa_W]\cup u_t=0$, then $\sum a_j(\eta_j)_h=0$. Therefore we get that
$$
\langle \Theta_{W,\bar W} u_t,u_t\rangle_t \leq 0.
$$
In particular, if $\kappa_j\cup u_t=0$ in cohomology for all $j$, then the full curvature operator acting on $u_t$ is negative. This will play an important role in section 9.

\medskip

In the proof of Theorem 6.1 we may assume that the base dimension $m$ is equal to 1. Indeed, this means that the formula for the curvature form implicit in Theorem 6.1 holds when restricted to any disk, and then it must hold on all of the base.
We start from the classical fact that if $u_t$ is a holomorphic section of an hermitian holomorphic vector bundle $E$, then
\be
\frac{\partial^2}{\partial t\partial\bar t} \|u_t\|^2_t=
-\langle \Theta^E_{\partial/\partial t,\partial/\partial \bar t} u_t, u_t\rangle_t 
+ \|D'_{\partial/\partial t} u_t\|^2_t.
\ee
Thus we need to compute $i\ddbar \|u_t\|^2_t$ when $u_t$ is a holomorphic section and $u_t$ is harmonic on fibers. We choose some representative $\tilde u$. By Lemma 2.6 we have
$$
\|u_t\|^2_t=(-1)^q c_n p_*(\tilde u\wedge\overline{\tilde u} e^{-\phi} ).
$$
Since the pushforward commutes with $\dbar$ we get
$$
\dbar  p_*(\tilde u\wedge\overline{\tilde u} e^{-\phi} )=
p_*(\dbar\tilde u\wedge\overline{\tilde u} e^{-\phi})+ (-1)^n
p_*(\tilde u\wedge\overline{\dphi\tilde u} e^{-\phi}).
$$
We claim that the first term vanishes. For this we use that $\dbar\tilde u=dt \wedge \eta +d\bar t \wedge \nu$. The term containing $\eta$ gives no contribution for bidegree reasons. For the second term we use that $\nu$ is $\dbar$-exact on fibers, which gives that the fiber integral
$$
\int_{X_t} \tilde u\wedge \bar\nu e^{-\phi}=0
$$
since $\tilde u=u_t$ on $X_t$ is harmonic on $X_t$. Thus
$$
\dbar  p_*(\tilde u\wedge\overline{\tilde u} e^{-\phi} )=
 (-1)^n
p_*(\tilde u\wedge\overline{\dphi\tilde u} e^{-\phi}).
$$
Taking $\partial$ we find
$$
\ddbar p_*(\tilde u\wedge\overline{\tilde u} e^{-\phi} )=(-1)^np_*(\dphi\tilde u\wedge\overline{\dphi\tilde u}e^{-\phi})+ p_*(\tilde u\wedge\overline{\dbar\dphi\tilde u} e^{-\phi})=:A+B.
$$
For the first term we recall that
$$
\dphi\tilde u=dt \wedge \mu + d\bar t \wedge \xi.
$$
Clearly the mixed terms containing $\mu\wedge dt\wedge\bar\xi\wedge dt$ vanish,
so
\be
A=[p_*(\mu\wedge\bar\mu e^{-\phi}) - p_*(\xi\wedge\bar\xi e^{-\phi})] dt\wedge d\bar t.
\ee
For the $B$-term we use $\dbar\dphi=-\dphi\dbar +\ddbar\phi$. Hence

$$
B= -p_*(\tilde u\wedge\overline{\dphi\dbar\tilde u} e^{-\phi})-
  p_*(\ddbar\phi\wedge \tilde u\wedge\overline{\tilde u}e^{-\phi}).
$$ 
    
    As we have seen
    $$
    p_*(\tilde u\wedge\overline{\dbar\tilde u}e^{-\phi})=0.
    $$
    Taking $\dbar$ of this we get
    $$
    p_*(\dbar\tilde u\wedge\overline{\dbar\tilde u}e^{-\phi})= (-1)^{n+1}
    p_*(\tilde u\wedge\overline{\dphi\dbar\tilde u}e^{-\phi}),
    $$
    so
    $$
    B= -  p_*(\ddbar\phi\wedge \tilde u\wedge\overline{\tilde u}e^{-\phi})+(-1)^n
    p_*(\dbar\tilde u\wedge\overline{\dbar\tilde u}e^{-\phi}).
    $$
    Now we use that
    $$
    \dbar\tilde u=dt \wedge \eta+d\bar t \wedge \nu.
    $$
    As before, the mixed terms give no contribution so
    $$
    (-1)^n p_*(\dbar\tilde u\wedge\overline{\dbar\tilde u}e^{-\phi})=
    p_*(\eta\wedge\bar \eta e^{-\phi}) -p_*(\nu\wedge \bar\nu e^{-\phi}).
    $$
    Putting all this together (and using that $\Omega=-i\ddbar\phi$) we finally get

    $$
    i\ddbar\|u_t\|^2_t=
    $$
    $$(-1)^q c_n
    [p_*(\mu\wedge\bar \mu e^{-\phi}) -p_*(\xi\wedge\bar\xi e^{-\phi})+
      p_*(\eta\wedge\bar \eta e^{-\phi}) - p_*(\nu\wedge\bar \nu e^{-\phi})] idt\wedge d\bar t +
    $$
    $$
    (-1)^q c_np_*(\Omega\wedge\tilde u\wedge\overline{\tilde u} e^{-\phi}).
    $$

    Up to this point the formula holds for any choice  of representative. Now we take $\tilde u=\hat u$ to be the vertical  representative. Then, by Proposition 5.3 and the comments after it, all forms $\eta,\mu,\xi$ and $\nu$ are primitive so by Lemma 2.6 the pushforward terms can be expressed as norms. Since $\mu$ and $\nu$ are of bidegree $(p,q)$ whereas $\xi$ is $(p+1,q-1)$
    and $\eta$ is $(p-1,q+1)$ we get
    \be
    i\ddbar\|u_t\|^2_t= \left(\|\mu\|_t^2 +\|\xi\|^2_t -\|\nu\|^2_t -\|\eta\|^2_t
      +c_n(-1)^q\int_{X_t} c(\Omega)u_t\wedge \bar u_t e^{-\phi} \right) idt\wedge d\bar t,
    \ee
    where we have also used Proposition 5.2 in the last term.

    Here we  decompose $\mu=\mu_h+\mu_\perp$, where $\mu_h$ is harmonic and $\mu_\perp$ is orthogonal to harmonic forms. Then $\mu_h=D'_{\partial/\partial t} u_t$, so comparing with the general curvature formula (6.3) we get
      \be\label{bo_big_formula}
      \langle\Theta_{\partial/\partial t,\partial/\partial \bar t}\, u_t,u_t\rangle_t=
      -\|\mu_\perp\|^2 -\|\xi\|_t^2 +\|\nu\|^2_t+\|\eta\|^2_t - c_n(-1)^q\int_{X_t} c(\Omega)u_t\wedge \bar u_t e^{-\phi}.
      \ee

This formula contains two positive  contributions to the curvature, one coming from the norm of $\nu$ and one coming from the norm of $\eta$. We shall now see that the first of these can be eliminated and the second can be improved by replacing $\eta$ by its harmonic part. 

For this we decompose  the forms $\eta$, $\xi$ and $\nu$, as we did with $\mu$, into one harmonic part and one part orthogonal to harmonic forms. Since $\nu$ is exact, $\nu=\nu_\perp$. We now use formulas (5.3) and (5.4), which clearly also 
holds for the non harmonic parts of $\eta,\xi,\mu$ and $\nu$. By Proposition 2.4, we can now rewrite formula (6.6).

First we note that since $\dbar\eta_\perp=0$ and $\dbar\nu_\perp=0$,
  they are not only orthogonal to harmonic forms, but also to all of the kernel of $\dphi$. In the same way, $\mu_\perp$ and $\xi_\perp$ are orthogonal to the kernel of $\dbar$. Therefore we can apply Proposition 2.4 with $g= \eta_\perp$ and $f=-\mu_\perp$, and to $g=\nu_\perp$ and $f=-\xi_\perp$. Finally we use that
  $$
  \langle (\Box''+1)^{-1}\xi_\perp,\xi_\perp\rangle_t +\|\xi_h\|^2=
  \langle (\Box''+1)^{-1}\xi,\xi\rangle_t.
  $$
  Inserting this in formula (6.6) we obtain Theorem 6.1.

  We finally record the counterpart of Theorem 6.1 for positively curved bundles. This can be proved in the same way as Theorem 6.1 (using the analog of Proposition 2.4 for positively curved metrics as commented in section 2), or by Serre duality, using Theorem 6.1.

\begin{thm}\label{curvaturepos} 

 Let $L\to \X$ be a line bundle with metric $e^{-\phi}$ where $i\ddbar\phi>0$ on fibers and define the hermitian structure on $\H^{n-q,q}$ using the K\"ahler forms $\omega_t=i\ddbar\phi|_{X_t}$ on fibers. Then
$$   
 \langle\Theta_{\partial/\partial t_j,\partial/\partial \bar t_k}\, u_t,u_t\rangle_t  =
 $$
 $$
 \langle (\Box''+1)^{-1}\eta_j, \eta_k\rangle_t +\langle (\Box'' +1)^{-1}\nu_k,\nu_j\rangle_t 
  -\langle(\xi_k)_h,(\xi_j)_h\rangle_t +c_n(-1)^q\int_{X_t} c_{j,k}(\Omega)u_t\wedge \bar u_t e^{-\phi},
  $$
where $\mu,\xi,\eta$ are defined in formulas (5.1) and (5.2), using  a vertical representative of $u_t$. 
\end{thm}
 This formula generalizes Theorem 1.2 in   \cite{Berndtsson2}, which deals with the case $p=n, q=0$. Then $\xi=0$ (since it is of bidegree $(p+1,q-1)$) and $\nu=0$ since $\nu$ is of bidegree $(n,0)$ and therefore must vanish identically if $u_t$ is a holomorphic section. This corresponds to the simplest case of Theorem 6.1, when $(p,q)=(0,n)$. Then the only  contribution with 'bad sign' in that theorem, coming from $\eta_h$, must vanish for bidegree reasons.

\section{Fiberwise flat metrics}

In this section we will consider the case when the bundle $L$ has a metric $e^{-\phi}$ with $i\ddbar\phi\leq 0$ on $\X$ and $i\ddbar\phi=0$ on each fiber $X_t$. We then have to assume given an auxiliary K\"ahler form $\Omega$ on $\X$. (When we assumed $i\ddbar\phi<0$ on fibers  we could take  $\Omega=-i\ddbar\phi$, or $\Omega=-i\ddbar\phi+p^*(\beta)$ where $\beta$ is a local K\"ahler form on the base.) We denote by $\omega_t$, or just $\omega$,  the restriction of $\Omega$ to $X_t$.

The main difference as compared to the previous case is that the cohomology is no longer necessarily primitive (with respect to $\omega$). Let us first pause to discuss the notion of primitive cohomology classes in $H^{p,q}(X, L)$, where $(X,\omega)$ is a compact K\"ahler manifold and $(L,e^{-\phi})$ is an Hermitian holomorphic line bundle. If $p+q=n-k$, we say that a class $[u]$ in
$H^{p,q}(X, L)$ is primitive if the class $[\omega^{k+1}\wedge u]$ vanishes. This means that $\omega^{k+1}\wedge u=\dbar v$ for some $L$-valued form $v$ of 
degree $n+k+1$. By the pointwise Lefschetz theorem, $v=\omega^{k+1}\wedge v'$ and it follows that
$$
\omega^{k+1}\wedge (u-\dbar v')=0
$$
so we could  equivalently have taken as our definition that the class $[u]$ has a representative that is pointwise primitive. Thus, there is a natural notion of primitive classes also for cohomology with values in a line bundle. However, in general, the Lefschetz decomposition on forms does not induce a Lefschetz decomposition on cohomology.

We now assume that $\ddbar\phi=0$ on $X$. Then it follows from the Kodaira-Nakano formula that $\Box'_\phi=\Box''_\phi=:\Box$. Moreover, the commutator
$[\Box,\omega\wedge]$ vanishes. Indeed, this is  well known when $L$ is trivial and $\phi=0$. Since the statement is  local it  also holds when
$\ddbar\phi=0$,
since we can always find a local trivialization with $\phi=0$ then. The first formula  implies that a harmonic form $u$ satisfies $\partial^\phi u=0$, and the second formula shows that if a class $[u]$ is primitive, then its harmonic representative is (pointwise) primitive. Indeed
$$
\omega^{k+1}\wedge u_h=\dbar v
$$
implies $\omega^{k+1}\wedge u_h=0$ since   $\omega^{k+1}\wedge u_h$ is harmonic. 
For later use we also point out that if a cohomology class $[u]$ is primitive for the K\"ahler form $\omega$, it is also primitive for any other K\"ahler form $\omega'\in [\omega]$. Indeed, $\omega'=\omega +\dbar v$ so
$$
(\omega')^{k+1}\wedge u =\omega^{k+1}\wedge u +\dbar(v'\wedge u),
$$
so $(\omega')^{k+1}\wedge u$ is $\dbar$-exact if $\omega^{k+1}\wedge u$ is $\dbar$-exact.

We then have a Lefschetz decomposition of cohomology classes
$$ 
H^{p,q}=P^{p,q} +\omega\wedge P^{p-1,q-1} +\omega^2 \wedge P^{p-2,q-2} +...
$$ 
where $P^{l,m}$ denotes the space of primitive classes, from the corresponding decomposition of harmonic forms.   This implies
$$
H^{p,q}=P^{p,q}+\omega\wedge H^{p-1,q-1}
$$
and in particular $h^{p,q}=p^{p,q}+h^{p-1,q-1}$ for the dimensions of the corresponding spaces.

We now apply this fiberwise  to our fiber space $p:\X\to\mathcal{B}$, where $h^{p,q}_t$ is constant over $\mathcal{B}$ and $p+q=n$. Then $h^{p-1,q-1}_t$ is upper semicontinuous and since $P^{p,q}$ is the kernel of a smoothly varying homomorphism its dimension is also upper semicontinuous. Therefore, since their sum is constant, both $p^{p,q}_t$ and $h^{p-1,q-1}_t$ are constant on $\mathcal{B}$. This implies in particular

\begin{prop}
  $\Pr^{p,q}$ is a smooth subbundle of $\H^{p,q}$.
\end{prop}
We also have

\begin{prop}
  $\Pr^{p,q}$ is a complex subbundle of $\H^{p,q}$.
  \end{prop}
\begin{proof}
  For this it is enough to show that if $u_t$ is in $\Gamma(\Pr^{p,q})$, then $D''u_t\in \Gamma^{0,1}(\Pr^{p,q})$. Let $\hat u$ be the vertical representative of $u_t$ with respect to the K\"ahler form $\Omega$, as defined in section 5. Recall that
  $$
  D'' u_t= \sum [\nu_j] d\bar t_j,
  $$
  where
  $$
  \dbar\hat u=\sum\eta_j\wedge dt_j+\sum \nu_j\wedge d\bar t_j.
  $$
  We have
  $$
  \nu_j\wedge dt\wedge d\bar t\wedge\Omega= \pm \dbar\hat u\wedge \widehat{d\bar t_j}\wedge dt\wedge \Omega=0,
  $$
  since $\Omega\wedge dt\wedge \hat u=0$. Hence $\nu_j$ is a  primitive representative of $[\nu_j]$.
\end{proof}
It follows in the same way that $\eta_j$ are primitive. 
As before we can now define a hermitian metric on our bundle $\Pr^{p,q}$ by
  $$
  \|[u]\|^2_t=(-1)^q c_n\int_{X_t} u_h\wedge\bar u_h e^{-\phi_t},
  $$
  i. e. the norm of a class is the norm of its harmonic representative, which is given by the above integral since the harmonic representative is primitive.

As we have seen,  if $[u_t]$ is a section of our bundle, then $\partial^\phi(u_t)_h=0$ on fibers. Hence, if $\hat u$ is the vertical representative of $(u_t)_h$,
  $$
  \partial^\phi\hat u=\sum\mu_j\wedge dt_j+\sum\xi_j\wedge d\bar t_j
  $$
  as before, and we see that $\mu_j$ and $\xi_j$ are primitive on fibers in the same way that we proved that $\nu_j|_{X_t}$ is primitive.  
  \begin{lma} $\ddbar\phi=p^*(C(\phi))$, where $C(\phi)$ is a $(1,1)$-form on the base.
  \end{lma}
  \begin{proof}
    Choose local coordinates $(t,z)$ on $\X$ such that $p(t,z)=t$. Since $\ddbar\phi$ vanishes on fibers we have
    $$
    \ddbar\phi=\ddbar_{t,\bar t}\phi +\ddbar_{t, \bar z}\phi+ \ddbar_{z, \bar t}\phi.
    $$
    Since $i\ddbar\phi\leq 0$ it follows from Cauchy's inequality that the mixed terms vanish, so
    $$
    \ddbar\phi=\sum \phi_{j k}(t,z)  dt_j\wedge d\bar t_k.
    $$
    Finally, the condition that $\ddbar\phi$ is $d$-closed gives that the coefficients are independent of $z$, which proves the lemma.
    \end{proof}
  \begin{prop}
    On each fiber
    $$
    \dbar\mu_j+\partial^\phi\eta_j=0
    $$
    and
    $$
    \dbar\xi_j+\partial^\phi\nu_j=0.
    $$
  \end{prop}
  \begin{proof}
    We have
    $$
    (\dbar\mu_j+\partial^\phi\eta_j)\wedge dt\wedge d\bar t=
    (\dbar\partial^\phi+\partial^\phi\dbar)\hat u\wedge \widehat{dt_j}\wedge d\bar t= \ddbar\phi\wedge \hat u\wedge \widehat{ dt_j}\wedge d\bar t=0,
    $$
    by the previous lemma.
  \end{proof}
  For a moment we now assume that the base is one dimensional and get exactly as in section 6 that formula (6.5) still holds,
  $$
  \langle \Theta u,u\rangle=
 ( -\|\mu\|^2-\|\xi\|^2+\|\nu\|^2+\|\eta\|^2)dt\wedge d\bar t -c_n(-1)^q p_*(C(\phi)\wedge \hat u\wedge\overline{\hat u} e^{-\phi}).
  $$
  We then decompose the  forms $\mu, \xi,\eta$ and $\nu$ into a harmonic part and one part that is orthogonal to harmonic forms.  The harmonic part of $\nu$ is zero since the section is holomorphic, and the harmonic part of $\mu$ is also zero if we assume that $D' u_t=0$ at the given point. We then apply Proposition 2.5, and conclude that
 $\|\mu_\perp\|=\|\eta_\perp\|$ and $\|\xi_\perp\|=\|\nu_\perp\|$. 
Hence,  the curvature formula becomes
  $$
   \langle \Theta u,u\rangle=
  -\|\xi_h\|^2+\|\eta_h\|^2 -c_n(-1)^q p_*(C(\phi)\wedge \hat u\wedge\overline{\hat u} e^{-\phi}).
  $$
  Since this holds for the restriction of $\Theta$ to any line in the base we finally get the curvature formula:
  \begin{thm} Let $L\to\X$ be a line bundle with hermitian metric $e^{-\phi}$ where $i\ddbar\phi\leq 0$ and $i\ddbar\phi=0$ on fibers. Then the curvature of the $L^2$ metric on $\Pr^{p,q}$, $p+q=n$ is
    $$
    \langle\Theta u,u\rangle_t=
    \sum \langle (\eta_j)_h,(\eta_k)_h\rangle_t dt_j\wedge d\bar t_k-
    \sum \langle (\xi_j)_h,(\xi_k)_h\rangle_t dt_j\wedge d\bar t_k +
    C(\phi)\|u\|^2_t.
    $$
  \end{thm}
We could now have continued and applied similar arguments to $\H^{p-i,q-i}$. Summing up the results we get that the same curvature formula holds for the entire bundle $\H^{p,q}$. 
Notice that when $L$ is trivial and $\phi=0$, this is a classical formula of Griffiths, \cite{GT}. Here, however we will be content with the subbundle $P^{p,q}$ since that is enough for our applications.

\section{Estimates of the holomorphic sectional curvature and hyperbolicity}

\noindent Let $p:X\to Y$ be a surjective holomorphic map between two compact K\"ahler manifolds. As in the introduction part, we denote by $Y_0:=Y\setminus \Sigma$ the set of regular values of $p$, so that if $X_0:=p^{-1}(Y_0)$ the restriction $p: X_0\to Y_0$ becomes a proper submersion. Our goal in this section is to analyze the properties of the base $Y_0$ induced by the variation of the complex structure of the fibers
$X_y$ of $p$ together with the semi-positivity properties of the canonical bundle of the said fibers. Under certain hypothesis we will construct a Finsler metric on the subset $\B\subset Y_0$ (cf. section 1) whose holomorphic sectional curvature
is bounded from above by a negative constant. To this end
we follow the same line of arguments as in the work by To-Yeung, \cite{To-Yeung} and Schumacher, \cite{Schumacher}. We repeat them here (with some modifications) to see how they adapt in our more general setting. The main conclusion is  that when the base is one dimensional we get a metric on the base with curvature bounded from above by a negative constant, provided that our metrics satisfy an integral bound which is automatic in the K\"ahler-Einstein case. We also take the opportunity to include the case of families of Calabi-Yau manifolds, as this does not seem to have appeared in the literature yet.
For simplicity of formulation we will assume that the base is one dimensional.

\subsection{The canonically polarized case}\label{CaPo}

We start with the case of a family of canonically polarized manifolds; $p:\X\to \B$, and we assume here that $\B$ is of dimension 1.  We let $L= -K_{\X/\B}$ and let $h=e^{-\phi}$ be a smooth metric on $L$ with $\Omega=-i\ddbar\phi$ strictly positive on each fiber and semipositive on the total space. By a result of Schumacher this holds if $\phi$ is a (normalized) potential of the K\"ahler-Einstein metric on each fiber, but for the moment we make no such assumption. Let $t$ be a local coordinate on the base $\B$, $V$ be the horizontal lift of $\partial/\partial t$, and $\kappa$ be the section of the bundle with fibers $Z^{0,1}(X_t, T^{1,0}(X_t))$, defined by $\kappa_t=\dbar_{X_t}V$. As in the last paragraph of section 5, we let $u^0$ be the canonical trivializing section of $\H^{n,0}$ (the bundle with fibers $H^{n,0}(X_t, -K_{X_t})\sim \C$), and then define $u^q$ inductively by $u^q=(\kappa\cup u^{q-1})_h$. As we have seen in section 5, $u^q$ is a holomorphic section of $\H^{n-q,q}$.

It follows from Theorem 6.1 that
\be
\langle i\Theta u^q,u^q\rangle \leq (-\|\xi_h\|^2 +\|\eta_h\|^2) (idt\wedge d\bar t)
\ee
since $\langle (1+\Box)^{-1}\xi,\xi\rangle\geq \|\xi_h\|^2$, and $c(\Omega)\geq 0$ since we have assumed that $\Omega\geq 0$. Recall that $\eta=\kappa\cup u^q$, $\xi=\bar\kappa\cup u^q$, and the subscript $h$ means that we have taken the harmonic part. We start with the following important observation by To-Yeung, \cite{To-Yeung}.
\begin{prop}
  $$
  \|\xi_h\|^2\geq \frac{\|u^q\|^4}{\|u^{q-1}\|^2}
  $$
  if $u^{q-1}\neq 0$.
\end{prop}
\begin{proof}
  We have
  $$
  \langle\xi_h, u^{q-1}\rangle=\langle\bar\kappa\cup u^q,u^{q-1}\rangle=
  \|u^q\|^2.
  $$
  Hence Cauchy's inequality gives
  $$
  \|u^q\|^2\leq\|\xi_h\|\|u^{q-1}\|
  $$
  which gives the claim.
\end{proof}

Next we introduce the notation
$$
\phi_q=\log\|u^q\|^2.
$$
Then
$$
i\ddbar\phi_q\geq -\langle i\Theta u^q,u^q\rangle/\|u^q\|^2,
$$
so formula (8.1) together with the proposition gives that
\be
i\ddbar\phi_q\geq (e^{\phi_q-\phi_{q-1}} -e^{\phi_{q+1}-\phi_q})(idt\wedge d\bar t),
\ee
since $\|\eta_h\|^2=e^{\phi_{q+1}}$.

A few comments are in order.  Notice that $\phi_q$ is only locally defined since it depends on the local coordinate $t$ in the base. Changing local coordinates we see that $e^{\phi_q}$ transforms as a metric on the $q$:th power of the tangent bundle of $\B$, and  $e^{\phi_q-\phi_{q-1}}$ defines a metric on the tangent bundle. In particular $e^{\phi_1}$ is the generalized Weil-Petersson metric on $\B$, and it is the genuine Weil-Petersson metric when the metric $\phi$ on $K_{\X/\B}$ is the (normalized) K\"ahler-Einstein potential on each fiber. Notice also that $u^q$ may be identically 0 on all fibers. If this happens for some $q$ we let $m$ be the maximal $q$ such that $u^q$ is not identically 0. Since $u^m$ is a holomorphic section of a vector bundle it can then only vanish on an analytic set. Hence $e^{\phi_m}$ defines a singular metric on the $m$:th power of the tangent bundle of $\B$. We will assume that the family is effectively parametrized, which means that $m$ is at least equal to 1.

Multiplying (8.2) by $q$ and summing from 1 to $m$ we get
$$
i\ddbar\sum_1^m q\phi_q\geq \sum_1^m (e^{\phi_q-\phi_{q-1}})(idt\wedge d\bar t).
$$
If $a_q>0$ and $\sum a_q=1$ we get by the convexity of the exponential function that the right hand side here is greater than
$$
\sum_1^m a_q e^{\phi_q-\phi_{q-1}}\geq e^{\sum a_q(\phi_q-\phi_{q-1})}.
$$

  Now  take $a_q=c(m+(m-1)+...q)$, with $c$ chosen so that $\sum_1^m a_q=1$. Then $a_{q-1}-a_q=c(q-1)$ so since $\sum_1^m a_q(\phi_q-\phi_{q-1})= a_m\phi_m+(a_{m-1}-a_{m})\phi_{m-1} +...-a_1\phi_0 $  we get $\sum a_q(\phi_q-\phi_{q-1})= c\sum q\phi_q- a_1\phi_0$. Hence
  $$
  i\ddbar\sum_1^n cq\phi_q\geq c e^{c\sum q\phi_q} e^{-a_1\phi_0} (idt\wedge d\bar t).
  $$
  Moreover, since $e^{\phi_q-\phi_{q-1}}$ is a metric on the tangent bundle of $\B$,
  $$
  e^{\sum a_q(\phi_q-\phi_{q-1})}= e^{c\sum q\phi_q} e^{-a_1\phi_0}
  $$
  is also a metric on the tangent bundle of $\B$, and so is $ e^{c\sum q\phi_q}=:e^\Phi$, since $\phi_0$ is a function. In conclusion, there is a metric with fundamental form  $e^\Phi idt\wedge d\bar t$,  on  $\B$ which satisfies
$$
  i\ddbar\Phi\geq c e^\Phi e^{-a_1\phi_0}idt\wedge d\bar t,
  $$
  and thus has
curvature bounded from above by
  $$
  -ce^{-a_1\phi_0},
  $$
  and so by a fixed negative constant if
  $$
  \phi_0=\log c_n\int_{X_t} u^0\wedge \bar u^0 e^{-\phi}
  $$
  is bounded from above. The last requirement is automatic if $\phi$ is a normalized K\"ahler-Einstein potential, which means that
  $$
  u^0\wedge \bar u^0 e^{-\phi}=(i\ddbar\phi)^n/n!
  $$
  on each fiber.

\subsection{The Calabi-Yau case}\label{CY} 
Let $\Omega>0$ be an arbitrary K\"ahler form on the total space $X$ of the fibration $p$. As before $t$ denotes a local coordinate on the base, and we let $V$ be the horizontal lift of $\partial/\partial t$ with respect to $\Omega$. The Kodaira-Spencer representative $\kappa$, and the holomorphic sections $u^q$ of $\H^{n-q,q}$ are defined as in the canonically polarized case.
\medskip
   
\noindent We will use the fact that the fibers $X_t$ are Calabi-Yau
as follows.
\begin{prop}\label{m-Ber} There exists a unique metric $e^{-\phi_{X/Y}}$ on the relative canonical bundle $K_{X/Y}$ such that $\displaystyle \ddbar\phi_t:= \ddbar\phi_{X/Y}|_{X_t}= 0$ and 
   \be\label{cy0}
 c_n\int_{X_t} u^0\wedge \bar u^0 e^{\phi_t}=1,
 \ee
for each $t\in Y_0$. This metric satisfies  $i\ddbar\phi_{X/Y}\geq 0$. 
\end{prop}

\begin{proof} We first remark that by \cite{Tosatti}, for any fixed fiber $X_{t_0}$ there is a positive integer $m$ such that $\displaystyle mK_{X_{t_0}}$ is trivial. 
  This implies that $\displaystyle h^0(mK_{X_t})$ is then equal to 1 for every $t$. Indeed, in the complement of a closed analytic subset of the base $Y$ this is true by general semicontinuity arguments, cf. \cite{BaSt}.
  At special points $\tau$ we have $\displaystyle h^0(mK_{X_\tau})\geq 1$. But, since  $c_1(X_\tau)= 0$, a holomorphic section can never vanish. Therefore $mK_{X_\tau}$ is trivial, so $h^0(mK_{X_\tau})= 1$ at $\tau$ as well. 

  In conclusion there exists a positive integer $m$ such that $\displaystyle h^0(mK_{X_t})= 1$ for all $t\in Y_0$. Moreover, the group $H^0(X_t, mK_{X_t})$ is generated by a nowhere-vanishing section  which extends locally near $t$.
  \smallskip
  
 \noindent We show next that the metric $e^{-\phi_{X/Y}}$ with the properties stated in
 \ref{m-Ber} is simply the $m$-Bergman metric cf. \cite{BP}.

 Let $t\in Y_0$ be a regular value of the map $p$, and let $x\in X_t$ be a
 point of the fiber $X_t$. If $(z_1,\dots, z_n)$ is a coordinate system on $X_t$ centered at $x$, then we obtain a coordinate system $(z_1,\dots, z_n, t)$ of the total space $X_0$ at $x$ by adding $t$. This induces in particular a trivialization of the relative canonical bundle with respect to which the expression of the
 $m$-Bergman metric becomes
\be\label{cy1}
e^{\phi_{X/Y}(z, t)}=\frac{|f(z, t)|^{2/m}}{\int_{X_t}|u|^{2/m}}\, ,
\ee 
where the notations are as follows: $u$ is any non-zero section of $mK_{X_t}$, and $|u|^{2/m}$ is the corresponding volume element on $X_t$. As we have seen, $u$ admits an extension on the $p$-inverse image of a small open set centered at $t$. Locally near $x$ we write $\displaystyle u= f \frac{\left(dz\wedge dt\right)^{\otimes m}}{\left(dt\right)^{\otimes m}}$ for some holomorphic function $f$. This identifies with
$f\left(dz\right)^{\otimes m}$, so we see that we have
\be\label{cy2}
c_nu^0\wedge \bar u^0 e^{\phi_t(z)}= \frac{|f(z, t)|^{2/m}}{\int_{X_t}|u|^{2/m}}
d\lambda(z)
\ee
from which the normalization condition \eqref{cy0} follows.

It follows from \cite{BP} that $\phi_{X/Y}$ is plurisubharmonic (this will also be a consequence of the remark below), and by \eqref{cy1}  it is moreover smooth, since $f$ is nonvanishing . It also follows from\eqref{cy1} that $i\ddbar\phi_{X/Y}|_{X_t}= 0$. Uniqueness follows since any other metric which is flat on fibers must have the form $\phi=\phi_{X/Y} +p^*(\chi)$ where $\chi$ is a function on the base. Condition \eqref{cy0} then implies that $\chi=0$, so Proposition \ref{m-Ber} is proved.\end{proof}

{\bf Remark:}
K-I Yoshikawa \cite{Yoshikawa} has noted that  the curvature form $i\ddbar\phi_{X/Y}$ is actually the pullback of the K\"ahler form $\omega_{WP}$ of the Weil-Petersson metric on the base. Apart from its interest in itself, this gives another proof of the semipositivity. Below we indicate how Yoshikawa's result can be obtained from our Theorem 7.5 (which we had not observed before learning about Yoshikawa's work).

The Weil-Petersson metric for a Calabi-Yau family $p:\X\to Y$ is defined as follows (cf. \cite{Tian} and \cite{CLWang} for the case when the fibers have trivial canonical bundle):

Let as before $\Omega$ be a K\"ahler form on the total space $\X$. Normalising 
we may assume that
$$
\int_{X_y} \Omega^n/n! =1
$$
for all fibers. Fix a point $y$ in the base. By Yau's theorem (\cite{Yau}) there is a unique K\"ahler form $\omega_y \in [\Omega|_{X_y}]$ such that $\omega_y^n/n!$  is the unique Ricci-flat volume element on $X_y$. If $[\kappa]\in H^{0,1}(X_y, T^{1,0}(X_y))$ we define (slightly abusively)
$$
\|[\kappa]\|^2_{y, WP}:=\int_{X_y} |\kappa_h|^2_{\omega_y} \omega_y^n/n!,
$$
where $\kappa_h$ is the harmonic representative of $[\kappa]$ for the metric $\omega_y$. (More correctly, if $\kappa=\dbar V$, where $V$ is a lift of a field $\partial_t$ on the base, this is the Weil-Petersson norm of $\partial_t$.) A priori this definition depends on the choice of $\Omega$, but we shall see that it does not.

First we claim that the Weil-Petersson norm as we have defined it coincides with the norm that we have used above, i. e.
\be\label{WP}
\|[\kappa]\|^2_{y, WP}=-c_n\int_{X_y} (\kappa\cup u^0)_h\wedge \overline{(\kappa\cup u^0)_h} e^{\phi_{X/Y}}.
\ee
To see this we note that
\be\label{Vol}
\omega_y^n/n!= c_n u^0\wedge\overline{u^0} e^{\phi_{X/Y}}
\ee
since both sides are Ricci-flat volume elements of total volume 1. Consider the map
$$
\tau:\E^{0,1}(X_y, T^{1,0}(X_y)) \to \E^{n-1,1}(X_y, -K_{X_y}),
$$
defined by $\tau(\kappa)=\kappa\cup u^0$. This map sends $\dbar$-closed forms to $\dbar$-closed forms and $\dbar$-exact forms to $\dbar$-exact forms, so it induces a map $\tilde\tau$ on the corresponding cohomology groups. One also verifies that it is an isometry in the sense that
$$
\int_{X_y} |\kappa|^2_{\omega_y} \omega_y^n/n!=\int_{X_y} (\kappa\cup u^0)\wedge *\overline{(\kappa\cup u^0)} e^{\phi_{X/Y}}.
$$
(This uses \eqref{Vol}.) Hence, $\tau$ sends harmonic forms to harmonic forms. By our discussion of primitivity at the beginning of section 7, $(\kappa\cup u^0)_h$ is primitive, so \eqref{WP} follows from the previous formula.

We now apply Theorem 7.5 with $(p,q)=(n,0)$. The bundle $\P^{n,0}=\H^{n,0}$ is trivial even metrically (by \eqref{cy0} the holomorphic section $u^0$ has norm constant equal to 1), so its curvature is 0. Moreover, the term involving $\xi$ vanishes for bidegree reasons. It follows that, if $\partial_t$ is a tangent vector in the base, $V$ is a lift of $\partial_t$ and $\kappa=\dbar V$,
then
$$
\|[\kappa]\|^2_{y, WP} =C(\phi_{X/Y})(\partial_t,\bar\partial_t),
$$
so $C(\phi_{X/Y})$ is the K\"ahler form of the Weil-Petersson metric, $\omega_{WP}$, i. e.
\be\label{WP100}
i\ddbar\phi_{X/Y}=p^*(\omega_{WP}).
\ee
  In particular, the Weil-Petersson metric does not depend on the choice of $\Omega$.

  \qed

\noindent The rest of the argument now goes as in  the canonically polarized case. We put
  $$
  \phi_q=\log \|u^q\|^2.
  $$
 Then, $\phi_0=0$ by construction, and for $q\geq 1$ we get from Theorem 7.10 that
  $$
  i\ddbar\phi_q\geq (e^{\phi_q-\phi_{q-1}}-e^{\phi_{q+1}-\phi_q})(idt\wedge d\bar t).
  $$
  Defining $\Phi=c\sum_1^n q\phi_q$ again, we find that
  $$
  i\ddbar\Phi\geq ce^{\Phi}i dt\wedge d\bar t,
  $$
  and that $e^\Phi$ defines a metric on the tangent bundle of $Y$ with curvature bounded from above by a strictly negative constant.

\medskip

{\bf Remark:} We refer the reader to the preprint \cite{Xu} by the third author where it is showed that the curvature formula in the fiberwise flat case can be used in a different way to produce a Hermitian metric on $Y$ of strictly negative {\it bisectional} curvature.
\qed

\subsection{Hyperbolicity} A direct consequence of the results in the subsections \ref{CaPo} and \ref{CY} respectively, is the following:

\begin{thm}\label{consec} We assume that the hypothesis of Theorem \ref{hypb1} are satisfied.
Then there exists a subset $\B\subset Y_0$ such that:
\begin{enumerate}
\smallskip
  
\item[{\rm (1)}] The complement $Y_0\setminus \B$
  is a closed analytic set, say $S$.
\smallskip
  
\item[{\rm (2)}] There exists a Finsler metric on $\B$ locally bounded from above at each point of $S$ and whose holomorphic sectional curvature is bounded from above by $-C$.
\end{enumerate}
\end{thm}

Let us first show how to use the above theorem to prove Theorem \ref{hypb1}.

\begin{proof}[Proof of Theorem \ref{hypb1}] The proof is based on ideas in 
\cite{2Schumacher},  Proposition 12. Let $f: \mathbb C\to Y\setminus \Sigma$ be an entire curve in $Y\setminus \Sigma:=Y_0$. By the above theorem, there exists a smooth strictly negatively curved Finsler metric on a Zariski open subset $\B:=Y_0\setminus S$ which extends to a singular metric on $Y_0$.

First let us assume that $f$ does not lie in $S$, thus the image of $f$ contains at least one point in $\B$. Then two cases can happen: $f$ is a constant or there exists $t\in \mathbb C$ such that $f(t) \in \B$ and $df(t)\neq 0$. Let us prove that the second case can never happen: Consider the pull back along $f$ of our Finsler metric on $Y_0$, the second case would give a singular metric, say $e^{\phi(t)}i dt\wedge d\bar t$, on $\mathbb C$ such that $\phi$ is smooth on an open set in $\mathbb C$ and
$$
\phi_{t\bar t} \geq c e^{\phi(t)}
$$  
in the sense of distribution on $\mathbb C$, where $c$ is a positive constant.  But the Ahlfors--Schwarz Lemma (see page 17 in \cite{JPhyp}) implies that
$$
e^{\phi(t)} \leq \frac2{c} \frac{R^{-2}}{(1-|t|^2/R^2)^2} \to 0, \ \ \text{as} \ R\to \infty,
$$ 
thus $\phi\equiv 0$ and we get a contradiction (since $\phi$ should be locally bounded from above).
\smallskip

\noindent If the entire curve $f$ is contained in $S$, we argue as follows.
We can assume that $S$ is the Zariski closure of the image of $f$ (as if this is not the case, then we simply replace
$S$ with the said Zariski closure). If $S$ is non-singular, then we are done
by the argument in the first part of this proof. If not, let
$$\pi: \widehat S\to S$$
be a desingularization of $S$, and let $\displaystyle \widehat p: \widehat X\to \widehat S$ be the family obtained by base change and desingularization of the total space.

\noindent The properties that $\widehat p$ inherits from the initial map $p$ are as follows.

\begin{enumerate}

\item[(a)] \emph{The generic fiber of $\widehat p$ is isomorphic to the generic fiber of $p$.} This is indeed clear, since $S$ is not contained in the singular locus of $p$. 

\item[(b)] \emph{The map $\widehat p$ has maximal variation.} Again this is immediate, because the desingularization map $\pi$ is generically isomorphic.

\item[(c)] \emph{The entire curve $f$ lifts to $\widehat p$.} This is the case since $S$ is the Zariski closure of the image of $f$.   
\end{enumerate}

\noindent The points (a)-(c) together with the arguments already invoked in the first part of our proof allow us to conclude.
\end{proof}

\textbf{Remark}: If $M$ is a complex manifold, the \emph{Kobayashi-Royden infinitesimal metric} is a Finsler metric on $T_M$ defined as
\be\label{hyp1}
k(x, v):= \inf\{\lambda> 0 : \exists \, \gamma: \bD\to M, \gamma(0)= x, \gamma^\prime(0)= 1/\lambda v\}.
\ee
We take $M:= Y_0$, where we recall that $Y_0:=Y\setminus \Sigma$ was the set of regular values of the
map $p$.
Then the above prove implies that
\be\label{hyp2}
\frac{k(x, v)}{|v|}\geq C_0 > 0
\ee
where the norm of the vector $v$ in \eqref{hyp2} is measured with respect to our Finsler metric on a dense subset $\B$ of $Y_0$. By a result of H. Royden, cf. \cite{Roy} the function $k$ is upper
semi-continuous on $T_{Y_0}$. Thus \eqref{hyp2} also holds true on 
$Y_0$.

\begin{proof}[Proof  of Theorem \ref{consec}]
  We define the set $\B\subset Y_0$ as in the introduction, such that the restriction of $\H^{n-i, i}|_{\B}$ is a vector bundle, whose fiber at $t$ is
  $H^{n-i,i}(X_t, L|_{X_t})$. The Finsler metric defined in \ref{CaPo} and \ref{CY} has the required
curvature property. The only thing to be checked is that the functions $\phi_q$ are locally bounded from above near each point of $Y_0\setminus \B$. In fact, each function $\phi_q$ is the norm of the $q^{\rm th}$
contraction of the tautological section $u_0$ with $\dbar V$, where $V$ is the
horizontal lift of $\frac{\partial}{\partial t}$. The norm of the harmonic
representative of a cohomology class is smaller than the norm of any other representative, and since our fibration $p$ is smooth on $Y_0$, the
boundedness statement in (2) follows.
\end{proof}

\section{Extension of the metric}

\noindent
In this section the set-up is as follows. We are given a surjective map $p: X\to Y$
between two smooth projective manifolds, such that the fibers of $p$ are connected (this last condition is not really necessary...). This will be referred to as \emph{algebraic fiber space}. Such a map $p$ will not be a  
submersion in general, so we cannot use directly the results obtained in the previous sections. In order to formulate our next results, we will recall next the notion of \emph{logarithmic tangent bundle} cf. \cite{Del}, which offers the right context to deal with the (eventual) singularities of $p$. 

\subsection{Notations, conventions and statements}
Let $(W_\alpha)_{\alpha\in J}$ be a finite set of non-singular 
hypersurfaces of $X$ which have transverse intersections. The logarithmic tangent bundle 
$T_X\langle W\rangle$ is the vector bundle whose local frame in a coordinate set $U$ is given by
\begin{equation}\label{mmp1}
z_1\frac{\partial}{\partial z_1},\dots, z_k\frac{\partial}{\partial z_k}, \frac{\partial}{\partial z_{k+1}},\dots, 
\frac{\partial}{\partial z_{n}}
\end{equation}
where $z_1,\dots, z_n$ are coordinates defined on $U$, such that $W_j\cap U= (z_j= 0)$
for $j=1,\dots k$. Hence we assume implicitly that only $k$ among the hypersurfaces $W_\alpha$ intersect the coordinate set $U$. We remark that the logarithmic tangent bundle is a subsheaf 
of $T_X$; its dual is the logarithmic cotangent bundle $\Omega_X\langle W\rangle$.
\smallskip

\noindent Throughout the current section we will observe the
following conventions.

\noindent $\bullet$ We denote by $\Sigma\subset Y$ the set of singular values of $p$. Let $\Delta\subset \Sigma$ be the codimension
one subset of $\Sigma$. We assume that the components of $\Delta$ are
smooth and have transverse intersections.
Note that the closure of the difference $\Sigma\setminus \Delta$ is a set of codimension at least two.

\noindent $\bullet$ Let $B\subset Y\setminus \Sigma$ the Zariski open subset of $Y$ for which the dimension of the relevant cohomology groups \eqref{higher50} is constant.
\smallskip

\noindent $\bullet$ We assume that there exists a
Zariski open subset $\Y_0\subset Y$ whose complement has codimension greater than two, such that if we denote $\X_0:= p^{-1}(\Y_0)$ then the following holds.

For any $x_0$ and $y_0$ in $\X_0$ and $\Y_0$ respectively, such that $p(x_0)= y_0$ we have local coordinates $(z_1,\dots z_{n+m})$ centered at $x_0$ and $(t_1,\dots t_m)$ centered at $y_0$ with respect to which the map $p$  
is given by 
\be\label{mmp2}
(z_1,\dots, z_{m+n})\to (z_{n+1},\dots, z_{n+m-1}, z_{n+m}^{b_{n+m}}\prod_{j=1}^qz_j^{b_j})
\ee
where $b_l\geq 1$ above are strict positive integers and $q\leq n$. Moreover, locally near $y_0$ the set $\Delta$ is given by the equation $t_m= 0$. We are referring to \cite{Naka} for the local expression
\eqref{mmp2}.
\smallskip

\begin{remark} Given an algebraic fiber space
$f:{\widetilde X}\to Y$, part of the properties above can be achieved
modulo the modification of the total space $\pi: X\to {\widetilde X}$ along
the inverse image of the discriminant of $f$. Also, if the discriminant $\Delta$ is not a simple normal crossing divisor, we consider a log-resolution of $(Y, \Delta)$,
together with the corresponding fibered product.
Hence up to such transformations we can assume without loss of generality that the bullets above hold true
for the maps $p$ we are considering.
\end{remark}
\medskip

\noindent We will denote by $T_Y\langle\Delta\rangle$ the logarithmic tangent bundle of $(Y, \Delta)$ described locally as in \eqref{mmp1}; let $\Omega_Y\langle\Delta\rangle$ be its dual. 
\medskip

\noindent We have the following easy and well-known statement, consequence of the bullets above.
\begin{lma} We have a natural morphism of vector bundles 
\be\label{mmp3}
p^\star\left(\Omega_Y\langle\Delta\rangle\right)\to \Omega_X\langle W\rangle
\ee
which is defined and injective on $\X_0$.
\end{lma}

\begin{proof}
  The verification is immediate: in the complement of the divisor $\Delta$ things are clear. Let $x_0$ and $y_0$ in $X$ and $Y$ respectively
  be two points as in (ii). We have the coordinates $(z_i)$ and $(t_j)$ with respect to which the map $p$ writes as in \eqref{mmp2}. Then we have
\be\label{mmp50'} 
p^\star(dt_j)= {dz_{n+j}}
\ee
for $j= 1,\dots, m-1$ and
\be\label{mmp51}
p^\star\left(\frac{dt_m}{t_m}\right)= b_{n+m}\frac{dz_{n+m}}{z_{n+m}}+ \sum_{i=1}^qb_{i}\frac{dz_{i}}{z_{i}}.
\ee
Thus the lemma is proved, given that the right hand side of \eqref{mmp51}
is a local section of $\Omega_X\langle W\rangle$.
\end{proof}
\noindent We note that the map \eqref{mmp3} is an
injection of sheaves on $X$, but in order to obtain an injection
of vector bundles, in general we have to restrict to $\X_0$ (as one sees by considering the blow-up of a point in ${\mathbb C}^2$).

\medskip

\noindent Let $\Omega_{X/Y}\langle W\rangle$ be the co-kernel of \eqref{mmp3}, so that we have the exact sequence
\be\label{mmp4}
0\to p^\star\left(\Omega_Y\langle\Delta\rangle\right)\to \Omega_X\langle W\rangle\to \Omega_{X/Y}\langle W\rangle\to 0
\ee
on the open set $\X_0\subset X$.

\noindent For further use, we will give next the expression of a local frame of
the bundle $\displaystyle \Omega_{X/Y}\langle W\rangle$ with respect to the coordinates $z$ and $t$ considered above. Let $U\subset \X_0$ be the open set, on
which the functions $z_j$ are defined. Then the local frame of $\displaystyle
\Omega_{X}\langle W\rangle|_U$ are given by
\be\label{mmp8}
\frac{dz_1}{z_1},\dots, \frac{dz_q}{z_q}, dz_{q+1},\dots, dz_{n+m-1},
\frac{dz_{n+m}}{z_{n+m}}.
\ee
Thus, the local frame of $\displaystyle \Omega_{X/Y}\langle W\rangle|_U$ is given by the symbols
\be\label{mmp9}
\frac{dz_1}{z_1},\dots, \frac{dz_q}{z_q}, dz_{q+1},\dots, dz_{n}, \frac{dz_{n+m}}{z_{n+m}}
\ee
modulo the relation
\be\label{mmp10}
b_{n+m}\frac{dz_{n+m}}{z_{n+m}}+ \sum_{j=1}^qb_j\frac{dz_{j}}{z_{j}}= 0.
\ee
The edge morphism corresponding to the direct image of the dual of \eqref{mmp3} gives the 
analogue of the Kodaira-Spencer map in logarithmic setting, as follows
\be\label{mmp5}
ks: T_Y\langle \Delta\rangle\to {\mathcal R}^1p_\star T_{X/Y}\langle W\rangle.
\ee
It turns out that we have 
\be\label{mmp6}
T_{X/Y}\langle W\rangle\simeq \Omega_{X/Y}^{n-1}\langle W\rangle\otimes K_{X/Y}^{-1}\otimes {\mathcal O}(p^\star(\Delta)- W),
\ee
on $\X_0$, hence we can rewrite \eqref{mmp5} as follows
\be\label{mmp5'}
ks: T_Y\langle \Delta\rangle\to {\mathcal R}^1p_\star \left(\Omega_{X/Y}^{n-1}\langle W\rangle\otimes L\right).
\ee
provided that the twisting bundle equals
\be\label{mmp10'}
L:= K_{X/Y}^{-1}\otimes {\mathcal O}\left(p^\star(\Delta)- W\right).
\ee
Hence, this fits perfectly into the framework developed in the previous sections of this article.
\smallskip

In general, given any bundle $L$ and an index $i\geq 0$ we have a map
\be\label{mmp6'}
\tau^i: \cR^{i}p_\star\left(\Omega^{n-i}_{X/Y}\langle W\rangle\otimes L\right)\to 
\cR^{i+1}p_\star\left(\Omega^{n-i-1}_{X/Y}\langle W\rangle\otimes L\right)\otimes \Omega_Y\langle\Delta\rangle
\ee
obtained by contraction with \eqref{mmp5}. 
A local section of $\cR^{i}p_\star\left(\Omega^{n-i}_{X/Y}\langle W\rangle\otimes L\right)$ is represented by 
a $(0,i)$-form with values in the $(n-i)$ exterior power of the
relative logarithmic cotangent bundle twisted with $L$.
We couple this form with the $(0,1)$-form with values in the logarithmic tangent bundle of $X$ corresponding to a local section of  
$T_Y\langle \Delta\rangle$ via \eqref{mmp5}. The result is a $(0,i+1)$-form with values in the $(n-i-1)$ exterior power of the logarithmic cotangent bundle twisted with $L$.
This is the map \eqref{mmp6'}; as we see, it is only defined in the complement of a 
set of codimension at least two in $Y$.
\begin{df}\label{sheaf}
We denote by $\K^i$ the kernel of \eqref{mmp6'}. It is a subsheaf of $\displaystyle \cR^{i}p_\star\left(\Omega^{n-i}_{X/Y}\langle W\rangle\otimes L\right)$. Also, let $\K^i_f$ be the quotient of $\K^i$ by its torsion subsheaf.
\end{df}
\noindent We can assume that the restriction $\K^i|_B$ is a subbundle of $\displaystyle \cR^{i}p_\star\left(\Omega^{n-i}_{X/Y}\langle W\rangle\otimes L\right)\big|_B$, by shrinking the set $B$.

\begin{remark}\label{noyau}
In view of our main formula (cf. Theorem \ref{curvature}), it is clear that the curvature of $\K^i|_B$ has better chances to be semi-negative than the full bundle $\displaystyle \cR^{i}p_\star\left(\Omega^{n-i}_{X/Y}\langle W\rangle\otimes L\right)$. Indeed, if the section $[u]$ in \ref{curvature} is a local holomorphic section of 
$\K^i|_B$, then the fiberwise projection of the $\eta_j$'s on the space of harmonic forms is identically zero. 
\end{remark}

\medskip

\noindent The sheaves $\K^i$ have been extensively studied in algebraic geometry, cf. \cite{Kang}, \cite{BruneB} and the references therein.
We will be concerned here with their 
differential geometric properties.
To this end, we formulate the requirements below concerning
Hermitian bundle $(L, h_L)$.
\smallskip

\noindent $\left({\mathcal H}_1\right)$ We have
$\displaystyle i\Theta_{h_L}(L)\leq 0$ on $X$
and moreover $\displaystyle i\Theta_{h_L}(L)|_{X_y}= 0$ for each $y$ in the complement of some Zariski closed set.
\smallskip

\noindent $\left({\mathcal H}_2\right)$ We have
$\displaystyle i\Theta_{h_L}(L)\leq 0$ on $X$
and moreover there exists a K\"ahler metric $\omega_Y$ on $Y$ such that we
have $\displaystyle i\Theta_{h_L}(L)\wedge p^\star\omega_Y^m\leq -\varepsilon_0\omega\wedge p^\star\omega_Y^m$ on $X$.
\smallskip

\noindent Thus the first condition requires $L$ to be semi-negative and
trivial on fibers, whereas in $\left({\mathcal H}_2\right)$ we assume that $L$
is uniformly strictly negative on fibers, in the sense that  we have   
\be\label{mp0301}
i\Theta_{h_L}(L)|_{X_y}\leq - \varepsilon_0\omega|_{X_y}
\ee
for any regular value $y$ of the map $p$.

\medskip

\noindent We recall that we denote by $\K^i_f$ the quotient of $\K^i$ by its torsion subsheaf (cf. Definition \ref{sheaf}). In this context we have the following result.

\begin{thm}\label{kernels, I}
  Let $p:X\to Y$ be an algebraic fiber space, and let $(L, h_L)$ be a Hermitian
  line bundle which satisfies one of the hypothesis
  $\left({\mathcal H}_i\right)$ above. We assume that the restriction of $h_L$
  to the generic fiber of $p$ is non-singular. Then: 
  \begin{enumerate}
  
  \item[(a)] For each $i\geq 1$ the 
  bundle $\displaystyle \K^i_f$ admits a semi-negatively curved singular Hermitian metric. Moreover, this metric is smooth on a Zariski open subset of $Y$.
  \smallskip
  
    \item[(b)] We assume that curvature form of $L$ is smaller than $-\varepsilon_0p^\star\omega_Y$ on the $p$-inverse image of some open subset $\Omega\subset Y$ of $Y$. Then the metric on $\displaystyle \K^i_f$ is strongly negatively curved on $\Omega$ (and semi-negatively curved in the complement of a codimension greater than two subset of $Y$).
\end{enumerate}    
  
 \end{thm} 

\smallskip

\noindent The method of proof of Theorem \ref{kernels, I} also gives the
following statement, which is potentially important in the analysis of families of holomorphic disks tangent to 
the pair $(Y, \Delta)$. Prior to stating this result, we define the ``iterated Kodaira-Spencer map'', cf. \cite{Schumacher}
\be\label{mmp58}
ks^{(i)}: \Sym^iT_Y\langle \Delta\rangle\otimes
{\mathcal R}^0p_\star\left(\Omega^{n}_{X/Y}\langle W\rangle\otimes L\right) \to {\mathcal R}^ip_\star\left(\Omega^{n-i}_{X/Y}\langle W\rangle\otimes L\right) 
\ee
by contracting successively the sections of $\displaystyle
\Omega^{\bullet}_{X/Y}\langle W\rangle\otimes L$
with the {$T_{X/Y}\langle W\rangle$-valued $(0,1)$-forms} given by $ks$ in
\eqref{mmp5} (we remark that the contraction operations are commutative, and this is the reason why the map \eqref{mmp58} is defined on $\Sym^i$ rather than $\otimes^i$ of the log tangent bundle of the base).
\medskip

\noindent In our next statement \ref{ks} we denote by $s_{\Delta}$ the section whose vanishing equals the support of $\Delta$. Also, we denote by $\Vert ks^{(i)}\Vert$
is the operator norm of the map \eqref{mmp58}.

\begin{prop}\label{ks} We assume that the singularities of the map $p$ are contained in the snc divisor $\Delta$. Let $\omega_{\P, X}$ and $\omega_{\P, Y}$ be two K\"ahler metrics with Poincar\'e singularities on $(X, W)$ and
  $(Y, \Delta)$ respectively. We assume that the metric $h_L$ of $L$ is smooth when restricted to the $p$-inverse image of a Zariski open subset of $Y$, and that its weights are bounded from below.
  Then for each $i\geq 1$ we have
\be\label{mmp52}
\Vert ks^{(i)}\Vert^{2/i}_y\leq C\log^N\frac{1}{|s_{\Delta}|_y^2},
\ee
where the
norm \eqref{mmp52} is induced by the metrics
$\omega_{\P, X}$ and $\omega_{\P, Y}$. The constants
$C, N$ are depending on everything but the point $y\in Y$. 
\end{prop}

\begin{remark} The proof of Proposition \ref{ks} is very similar to the arguments we will give for Theorem \ref{kernels, I}, and we have decided that we can afford to skip it. We will however highlight the main points after completing the arguments for \ref{kernels, I}.
\end{remark}  
\medskip

\noindent Assume that the bundle $(L, h_L)$ verifies one of the hypothesis
$\left({\mathcal H}_i\right)$, and that moreover we have
\be\label{mmp30}
H^0\left(X, \Omega^{n}_{X/Y}\langle W\rangle\otimes L\right)\neq 0.
\ee
If $\sigma$ is a holomorphic section of $\displaystyle \Omega^{n}_{X/Y}\langle W\rangle\otimes L$, then we have a holomorphic map
\be\label{mmp31} 
ks^{(1)}_\sigma: \O_Y\to 
\cR^{1}p_\star\left(\Omega^{n-1}_{X/Y}\langle W\rangle\otimes L\right)\otimes \Omega_Y\langle\Delta\rangle.
\ee
induced by the map $ks^{(1)}$ cf. \eqref{mmp58}. 

\medskip

\noindent We formulate yet another hypothesis.
\smallskip

\noindent $\left({\mathcal H}_3\right)$ There exists a holomorphic
section $\sigma$ of $\displaystyle \Omega^{n}_{X/Y}\langle W\rangle\otimes L$
such that the map $\displaystyle ks^{(1)}_\sigma$
is not identically zero when restricted to a
non-empty open subset of $\B$.
\medskip

\noindent The following result is a corollary of Theorem \ref{kernels, I}.

\begin{thm}\label{kernels, II}
Let $p:X\to Y$ be an algebraic fiber space, and let $(L, h_L)$ be a Hermitian
  line bundle which satisfies one of the hypothesis
  $\left({\mathcal H}_j\right)$ for $j=1, 2$. 
Moreover, we assume that $\left({\mathcal H}_3\right)$ is equally satisfied.
Then there exists $s\leq n= \dim (X_y)$ and a non-trivial map
\be\label{mmp31'} 
\K^{s\star}_f\to \Sym^s \Omega_Y\langle \Delta\rangle.
\ee
In addition, if the curvature of $L$ is smaller than
$-\varepsilon_0p^\star\omega_Y$ on the pre-image of some non-empty 
open subset $V$ of $Y$,
then there exists an ample line bundle $A$ on $Y$ such that the bundle $$\otimes^M\Omega_Y\langle  \Delta\rangle\otimes A^{-1}$$ has a (non-identically zero) global section, for some $M\gg 0$. 
\end{thm}
\medskip

\noindent Our next results will show that in some interesting geometric circumstances a version of the 
hypothesis $\left({\mathcal H}_2\right)$ is verified. 

If the canonical bundle of the generic fiber of $p: X\to Y$ is ample, then we say that the family defined by $p$ is \emph{canonically polarized}. If the
Chern class of the canonical bundle of the generic fiber of $p$ equals zero,
then we say that $p$ defines a \emph{Calabi-Yau} family. 
Next, a canonically polarized or Calabi-Yau family $p$ has \emph{maximal variation} if the Kodaira-Spencer
map \eqref{mmp5} is injective when restricted to a non-empty open subset of
$Y_0$. 
\medskip

\noindent Then we have the following result, which provides a (geometric) sufficient condition under which
the anti-canonical bundle of $p$ has strong curvature properties. 

\begin{thm}\label{kernels}
Let $p:X\to Y$ be a family of canonically polarized manifolds of maximal variation, and let $\omega$ be a reference metric on $X$. Then there exists an effective divisor $\Xi$ in $X$ and a metric $h_{X/Y}$ on the twisted relative canonical bundle 
$K_{X/Y}+ \Xi$ such that the following hold true: 
\begin{enumerate}
  
\item[(i)] The curvature corresponding to $h_{X/Y}$ is greater than $\ep_0\omega$ for some $\ep_0> 0$ and
  the restriction $h_{X/Y}|_{X_y}$ is non-singular for all $y$ in the complement of a Zariski closed subset of $Y$.
  \smallskip

\item[(ii)] The codimension of the direct image $p(\Xi)$ is greater than two.
\end{enumerate}  
\end{thm}

\begin{remark}Given that $\displaystyle K_{X_y}$ is ample, a natural choice for
  the construction of a positively curved metric on $K_{X/Y}$
  would be the fiberwise KE metric, cf. previous section --this is well defined in the
  complement of the singular locus of $p$, and it extends across the singularities, cf. \cite{Paun}. However, it is not clear for us whether it is possible to obtain a \emph{useful} lower bound of the eigenvalues
  of the KE metric with respect to a fixed K\"ahler metric on $X$ as we are
  approaching the singular locus of $p$. And this information is critical in the process of extending the metric, as we will see below. 
\end{remark}

\medskip

\noindent In order to derive an interesting consequence of Theorem \ref{kernels} 
we recall the following result for which we refer to \cite{CP}, Theorem 2.3.
 
\begin{thm}\label{bo_powers}\cite{CP}
We assume that the hypothesis of Theorem \ref{kernels} are satisfied. Then there exists a metric $h_{X/Y}^{(1)}$
on the relative canonical bundle $K_{X/Y}$ with the following properties.
\begin{enumerate}

\item[(i)] The restriction of $h_{X/Y}^{(1)}$ to the fiber $X_y:= p^{-1}(y)$ is smooth. It is induced by the sections of $\displaystyle mK_{X_y}$ for some $m\gg 0$, for all $y$ in the complement of a proper algebraic subset $Z\subset Y$. 
\smallskip

\item[(ii)] The corresponding curvature current $\Theta$ is positive definite on each compact subset of $p^{-1}(Y\setminus Z)$, and moreover we have
\be\label{031001mp}
\Theta\geq \sum(t^j- 1)[W_j]
\ee
in the sense of currents on $X$.
\end{enumerate}
\end{thm}

\noindent In the statement \eqref{bo_powers} we denote by $t^j$ the multiplicity of the hypersurface $W_j$ in the inverse image $p^{-1}(\Delta)$.  By combining Theorem \ref{kernels} and Theorem \ref{bo_powers}, we obtain the following result.

\begin{cor}\label{approx}
Let $\eta> 0$ be a positive real number. We assume that the hypothesis of Theorem \ref{kernels}
are satisfied, and we define the following family of metrics:
\be\label{031002mp}
\psi_{X/Y}^{(\eta)}:= (1-\eta)\varphi_{X/Y}+ \eta\varphi_{X/Y}^{(1)}- \sum(t^j-1)\log |f_j|^2
\ee
on the bundle $K_{X/Y}- \sum (t^j-1)W_j$, where $\varphi_{X/Y}$ and $\varphi_{X/Y}^{(1)}$ are the weights given by Theorem \ref{kernels} and 
Theorem \ref{bo_powers}, respectively, and $f_j$ is a local equation of the hypersurface $W_j$. 
The resulting metrics $h^{(\eta)}$ have the following properties.
\begin{enumerate}

\item[(a)] The curvature current $\Theta_\eta$ verifies 
$$\Theta_\eta\geq (1-\eta)\ep_0\omega- (1-\eta)\sum(t^j- 1)[W_j]- \eta[\Xi].$$
\smallskip

\item[(b)] Each of the metrics $h^{(\eta)}$ is smooth when restricted to the generic fiber of $p$.
\smallskip

\item[(c)] For each compact subset $K\subset X\setminus p^{-1}(\Delta)$ we have $\Theta_\eta|_K\geq C_K\omega|_K$, where $C_K> 0$ is independent of $\eta$.

\end{enumerate}
\end{cor}

\noindent Thus, in the set-up of Theorem \ref{kernels} the point (a) of Corollary \ref{approx} shows that a \emph{version} of hypothesis
$\left({\mathcal H}_2\right)$ is satisfied. Indeed, we consider the bundle
$L:= -K_{X/Y}+ \sum (t^j-1)W_j$, endowed with the metric induced by
$h^{(1)}$. This metric does not necessarily verifies the requirements in $\left({\mathcal H}_2\right)$. Nevertheless, it is the limit of $h^{(\eta)}$
as $\eta\to 1$, and the restriction to fibers of $\Theta_\eta$ is uniformly positive, cf. (a).     
We will show at the end of this section that that this suffices to infer that Theorem 
\ref{kernels, I} holds true. 

\medskip

\noindent The following important result due to Viehweg-Zuo (cf. \cite{VZ03}) can now be seen as a direct consequence of the results \ref{kernels, I},
\ref{kernels, II} and \ref{kernels}.

\begin{thm}\label{VZ}\cite{VZ03}
  Let $p:X\to Y$ be a family of canonically polarized manifolds of maximal variation. Then there exists a positive integer $q\leq \dim(Y)$ such that the bundle $\Sym^q\Omega_Y\langle \Delta\rangle$ contains a non-trivial big coherent subsheaf.  
\end{thm}
\noindent We refer to \emph{loc. cit.} for the notion of ``big subsheaf''
$\mathcal F$: it implies that for any ample line bundle $H$ the bi-dual of $\Sym^m\mathcal F\otimes H^{-1}$ is generated by global sections on some open subset $U\subset X$ for all $m\gg 0$.
\medskip

\noindent As a by-product of our methods it turns out that a completely similar result holds in the context of Calabi-Yau families.

\begin{thm}\label{CY17}
  Let $p:X\to Y$ be a Calabi-Yau family. Assume that $p$ has maximal variation. Then there exists a positive integer $q\leq \dim(Y)$ such that the bundle $\Sym^q\Omega_Y\langle \Delta\rangle$ contains a non-trivial big coherent subsheaf.  
\end{thm}

\medskip

\medskip

\noindent In the following subsections we will establish the 
results stated above.

\subsection{Proof of Theorem \ref{kernels, I}}
Our first task will be to construct a metric on the quotient
\be\label{corrmp1}
\K^i_f:= \K^i/T(\K^i)
\ee
of the kernel by its torsion subsheaf. This will be naturally
induced by the metric on the direct image, thanks to the following
simple observation.

\begin{lma}\label{torsquot} The support
  of the torsion subsheaf $T(\K^i)$ is contained in $Y\setminus B$.
\end{lma}
\noindent This is clear, since $\K^i|_B$ is a vector bundle.


\smallskip

\noindent Next, the sheaf $\K^i_f$ is coherent and
torsion-free, therefore it is locally free on an
open subset $\Y_0$ whose codimension in $Y$ is at least two. We define next a metric on the restriction $\K^i_f|_B$ as follows.

Let $s_1, s_2$ be local sections of $\K^i_f$ defined on an open set $V$ such that the sequence
\be\label{corrmp3}
0\to \Gamma\left(V, T(\K^i)\right)\to\Gamma\left(V, \K^i\right)\to
\Gamma\left(V, \K^i_f\right)\to 0
\ee
is exact. We consider two sections $u_1, u_2\in \Gamma\left(V, \K^i\right)$ projecting into
 $s_1$ and $s_2$, respectively via \eqref{corrmp3}. Then for each
$t\in B\cap V$ we define
\be\label{corrmp4}
\langle s_1, s_2\rangle_t:= \langle u_1, u_2\rangle_t. 
\ee
We note that this does not depends on the $u_i$'s by Lemma \ref{torsquot}, and therefore we have a well-defined Hermitian structure on $\K^i_f|_B$.
\smallskip

\noindent Next, we will use
the Leray isomorphism in order to construct a representative of a local section $[u]$ of $\K^i$. Since this is slightly different from the convention adopted in the previous sections, we give a few precisons in what follows.

Assume that $[u]$ is defined on a
small enough open subset $V\subset Y$ containing a smooth point of the divisor $\Delta$. In particular, $[u]$ corresponds to an element
of the cohomology group
\be\label{corrmp10}
H^i\left(p^{-1}(V), \Omega^{n-i}_{X/Y}\langle W\rangle\otimes L|_{p^{-1}(V)}\right).
\ee
Let $\U= (U_\alpha)$ be a finite cover of $X$, such that the intersections $U_{I}:= U_{\alpha_0}\cap\dots\cap U_{\alpha_s}$ are
contractible for all multi-indexes $I= (\alpha_0<\dots< \alpha_s)$ having $s+1$ components and all $s$.

The section $[u]$ corresponds to a collection of holomorphic sections $(u_I)$ of the bundle
$$\displaystyle \Omega^{n-i}_{X/Y}\langle W\rangle\otimes L|_{U_I}$$ whose Cech co-boundary is equal to zero (i.e. it is a cocycle). Here the multi-index $I$
has $i+1$ components. We can assume that the map
\be\label{corrmp11}
\Omega^{n-i}_{X}\langle W\rangle\otimes L|_{U_I}\to \Omega^{n-i}_{X/Y}\langle W\rangle\otimes L|_{U_I}
\ee
is surjective at the level of sections.
Let $\wt u_I$ be a lifting of the section $u_I$ via the map
\eqref{corrmp11}. Each $\wt u_I$ is an $L$-valued $n-i$-form with
log poles on $U_I$. 
We note that in general the collection $(\wt u_I)$ is not necessarily a cocycle. However, this is the case for the
restriction $\displaystyle (\wt u_I|_{U_I\cap X_t})$ for each $t\in B$,
and this identifies with the class $[u_t]$.

We apply the Leray isomorphism procedure to $(\wt u_I)$; this gives a
$(0, i)$-form with values in $\displaystyle \Omega^{n-i}_{X}\langle W\rangle\otimes L|_{p^{-1}(V)}$ which we denote by $\wt u$. The exact formula is as follows
\be\label{corrmp12}
\wt u= \sum_I\rho_{\alpha_s}\wt u_{\alpha_0\dots \alpha_s}\dbar \rho_{\alpha_0}\wedge \dots\wedge \dbar \rho_{\alpha_{s-1}}
\ee
where $(\rho_\alpha)$ is a partition of unit subordinate to the Leray cover $\U$. In the usual formula, we do not have the factor $\rho_{\alpha_s}$ in \eqref{corrmp12}, since the forms defined on 
the open sets $U_{\alpha_s}$ coincide on intersections. As already mentioned, here $(\wt u_I)$ is not necessarily a cocycle, but this is the case for the projection $(u_I)$. Since $\sum \rho_j=1$, the image of $\wt u$ on
the space of $(0,i)$-forms with values in $\Omega^{n-i}_{X/Y}\langle W\rangle\otimes L|_{p^{-1}(V)}$ is $\dbar$-closed and it represents the class $[u]$.
\smallskip

\noindent -- \emph{Convention.}-- For the rest of this section, we will call $\wt u$ a representative of $[u]$. 
Also, the representative of a section $u$ of $\K^i_f$ will be by definition the representative of a section of $\K^i$ projecting into
$u$.
\smallskip

\noindent Let $\wt u$ be a fixed representative of a local
section of the bundle $\K^i_f$, which is defined near a smooth point $y_0\in \Delta$. Thus $\wt u$ is a $(0, i)$-form with values in the bundle
$\displaystyle \Omega^{n-i}_{X}\langle W\rangle\otimes L|_{p^{-1}(V)}$,
such that the following properties are satisfied for any $t\in B\cap V$.

\noindent $\bullet$ 
The restriction of the form $\wt u$ to the fibers $X_t$ 
is $\db$-closed. 
\smallskip

\noindent $\bullet$ The cup-product of $\wt u$ with any element in the image of $ks$ is cohomologically zero
when restricted to $X_t$.
\smallskip

\noindent We first
assume that the bundle $(L, h_L)$ verifies the hypothesis $\left({\mathcal H}_2\right)$.
 Let
$\theta:= -\Theta_{h_L}(L)$ be the corresponding curvature form. Then we have
\be\label{mmp11}
\sqrt{-1}c_{jk}(\theta)dt_j\wedge d\overline{t}_k\geq 0
\ee
and given that $\wt u$ is a section of the kernel, the fiberwise projection of the (troublesome) forms $\eta_j$ on the space of harmonic forms is equal to zero (as above, we assume that $t\in B$).
Thus, for any local holomorphic section $u$ of $\K^i_f$ defined on an open set $V$ centered at the point $y_0$ the function 
\be\label{mmp12}
t\to \log\Vert [u]\Vert_t^2
\ee
is psh on $V\setminus \Delta$, as a consequence of Theorem \ref{curvature} (the bracket notation in \eqref{mmp12} has the same meaning as in section 2). We will show next that we have
\be\label{mmp13}
\sup_{t\in V\setminus \Delta}\log\Vert [u]\Vert_t^2<\infty
\ee
and then Theorem \ref{kernels, I} follows as a consequence of elementary
properties of psh functions.

In order to establish \eqref{mmp13}, by the definition of the norm of a section of
${\mathcal H}^{p, q}$ we have the inequality
\be\label{mmp14}
\Vert [u]\Vert_t^2\leq \int_{X_t}\left|\wt u_{|X_t}\right|^2_{\theta, h_L}\theta^n
\ee
for each $t\in V\setminus \Delta$. Hence it would be enough to bound the right hand side term of \eqref{mmp14}.
\smallskip

\noindent Unfortunately we cannot do this directly (as the following computations will show, the technical reason is that we do not have an upper bound for the eigenvalues of $\theta$ at our disposal) and we will
proceed as follows.

 By hypothesis $\left({\mathcal H}_2\right)$ the form $\theta$ is
semi-positive on $X$ and greater than $\varepsilon_0\omega$ on the fibers of $p$. For each $j$ let $s_j$ be the global section whose set of zeroes is the hypersurface $W_j$ (recall that the $W_j$ are the support of the inverse image $p^{-1}(\Delta)$). They correspond locally on the open set $U$ (cf. \eqref{mmp2}) to the coordinates
$z_1,\dots, z_q, z_{m+n}$. We have the decomposition $\displaystyle U:= \bigcup U_j$, where $U_j\subset U$ corresponds to the set of points for which $|z_j|$ is minimum among $|z_1|,\dots, |z_q|, |z_{n+m}|$.
If $j=1,\dots, q$ then on the open set $U_j\cap X_t$ of the fiber $X_t$ we take the local coordinates
$z_1,\dots z_{j-1}, z_{n+m}, z_{j+1},\dots, z_n$. On the set $X_t\cap U_{n+m}$ the coordinates are 
$z_1,\dots, z_n$.

\noindent There exists a constant $C> 0$ such that for each $\varepsilon> 0$, on $X_t\cap U_{n+m}$ we have
\be\label{mmp15}
\theta-\ep\sum_{j}dd^c\log\log\frac{1}{|s_j|^2}\geq
C\sqrt{-1}\left(\sum_{j=1}^ndz_j\wedge d\ol z_j+ \ep\sum_{j=1}^q\frac{dz_j\wedge d\ol z_j}
{|z_j|^2\log^2|z_j|^2}\right).
\ee
\smallskip

\noindent The inequality \eqref{mmp15}
is a direct consequence of the fact that $\theta$ is greater than
$\varepsilon_0\omega$ on the fibers of $p$, by hypothesis. 

In order to apply Theorem \ref{curvature}, we perturb the metric of $L$ as follows
\be\label{mmp16}
h_{L,\ep}:= \left(\prod_j \log^\ep\frac{1}{|s_j|^2}\right)^{-1}
e^{C\sqrt{\ep}|t|^2}h_L
\ee
and let $\theta_\ep$ be the opposite of resulting curvature form (whose restriction to fibers is none other than
the left hand side of \eqref{mmp15}). The constant $C\gg 0$
in \eqref{mmp16} is chosen so that $\theta_\ep$ still verifies $\left({\mathcal H}_2\right)$, for each positive $\ep> 0$. This is quickly seen as follows. We have
\be\label{mmp2016}
-dd^c\log\log\frac{1}{|s_j|^2}= \frac{\sqrt{-1}}{\pi}\frac{\langle D's_j, D's_j\rangle}{|s_j|^2
\log^2|s_j|^2}+ \frac{1}{\log|s_j|^2}\Theta(W_j),
\ee
so we remark that the negativity induced by the first factor of 
\eqref{mmp16} in the expression of $\theta_\ep$ is of order $\displaystyle \sum_j\frac{\ep}{\log1/|s_j|^2}\omega$. This is the reason why we introduce the weight $C\sqrt{\ep}|t|^2$ in \eqref{mmp16}: given the strict positivity of $\theta$ on the fibers of $p$, the curvature of $\displaystyle h_{L, \ep}$ will be semi-negative on $X$.

\noindent In conclusion, we can change the metric of $L$ in order to have Poincar\'e singularities on the divisor $W= \sum W_j$.

\smallskip

\noindent We clearly have

\be\label{mmp17}
\lim_{\ep\to 0}\Vert [u]\Vert^2_{\ep, t}= \Vert [u]\Vert^2_t
\ee
for every $t\in V\setminus \Delta$.  
In \eqref{mmp17} above we denote by $\displaystyle \Vert \cdot\Vert^2_{\ep, t}$ the norm induced by the 
perturbed metric of $L$. The equality \eqref{mmp17} is a consequence of the usual elliptic theory, cf. \cite{KSp}, applied for each $t$.
\medskip

\noindent We formulate the following claim,
where we recall that $\wt u$ is a representative of the section $u$.
\smallskip

\noindent{\bf Claim.} For each $\ep> 0$ there exists a constant $C_\ep> 0$
such that we have
\be\label{mmp18}
\sup_{t\in V\setminus \Delta}|t_m|^{2\ep}\int_{X_t}\left|\wt u_{|X_t}\right|^2_{\theta_\ep, h_{L, \ep}}\theta^n_\ep
\leq C_\ep< \infty.
\ee
\medskip

\noindent If we are able to show that \eqref{mmp18} holds true, then we are done (despite of the fact that the constant $C_\ep$ above may not be bounded as $\ep\to 0$). Indeed, this would imply that the log of the expression 
\be\label{mmp22}
t\to |t_m|^{2\ep}\Vert [u]\Vert^2_{\ep, t}
\ee
defines a psh function on $V$ for each $\ep > 0$. We write the corresponding mean inequality 
\be\label{mp0309}
\log \left(|\tau|^{2\ep}\Vert [u]\Vert^2_{\ep, {(t', \tau)}}\right)\leq \ep C+ 
\int_0^{2\pi} \!\!\!\log \Vert [u]\Vert^2_{\ep, {(t', \tau_\theta)}}\frac{d\theta}{2\pi}
\ee
where we use the notation $\tau_\theta:= \tau+ re^{\sqrt{i}\theta}$ under the integral sign in \eqref{mp0309}, and
$0< r\ll 1$.

Now, the concavity of the log function implies
\be\label{mpdec1}
\int_0^{2\pi} \!\!\!\log \Vert [u]\Vert^2_{\ep, {(t', \tau_\theta)}}\frac{d\theta}{2\pi}\leq \log\left(\int_0^{2\pi} \!\!\!\Vert [u]\Vert^2_{\ep, {(t', \tau_\theta)}}\frac{d\theta}{2\pi}\right),
\ee
and we obviously have
\be\label{mpdec2}
\Vert [u]\Vert^2_{\ep, {(t', \tau_\theta)}}\leq \int_{X_t}|\wt u_{|X_t}|^2_{\theta_\ep, h_{L, \ep}}\theta^n_\ep 
\ee
where $t$ in \eqref{mpdec2} is equal to $\displaystyle (t', \tau_\theta)$. The \eqref{mpdec2} is true because the norm on the left hand side is given by the \emph{harmonic} representative and this is clearly smaller than the norm of the representative $\wt u$.

\noindent On the other hand, if $|t_m|= \delta_0> 0$, then we have
\be\label{mmp19}
\sup_{|t_m|= \delta_0,\hskip 2pt t\in V}\int_{X_t}|\wt u_{|X_t}|^2_{\theta_\ep, h_{L, \ep}}\theta^n_\ep
\leq C_0
\ee
\emph{uniformly with respect to $\ep$}, since we are ``away'' from the singularities of the map $p$. The mean inequality \eqref{mp0309}, combined with \eqref{mmp19}, \eqref{mmp17} give the uniform boundedness of the initial norm, which is what we wanted.
\smallskip

\noindent In conclusion, the argument is that we will be using the non-effective estimate \eqref{mmp18}
in order to infer that the function \eqref{mmp22} is log-psh. Then the uniform boundedness of our initial
norm follows by the mean inequality, together with the convergence as $\ep\to 0$ property \eqref{mmp17}.

\begin{proof} (of the Claim) To establish the claim 
we will use the local frame \eqref{mmp9} in order to do a few local computations. 
Without loss of generality we can assume that $z$ belongs to the coordinate set $U_{n+m}\subset U$. For each $t\in V\setminus \Delta$
the local coordinates on $X_t\cap U_{n+m}$ will be $z_1,\dots, z_n$.

Our representative $\wt u$ can be expressed as follows

\be\label{rev1}
\wt u= \sum_{I, J}\tau_{I\ol J}(z)\frac{dz_\alpha}{z_\alpha}\wedge dz_\beta\wedge d\ol z_J\otimes e_L
\ee
where $\alpha, \beta$ are multi-indexes whose union equals $I$, such that we have
$\alpha\subset \{1,\dots q, m+n\}$ and $\beta\subset \{q+1,\dots, m+n-1\}$. The set $J$ is contained in $\{1,\dots, m+n\}$ and it has length $i$. The coefficients
$\tau_{I\ol J}$ are differentiable functions, in particular uniformly bounded.

The restriction of $\wt u$ to $X_t\cap U_{n+m}$ reads as follows
\be\label{mmp23}
\wt u|_{X_t}= \sum_{I, J}\xi_{I\ol J}(z, t)\frac{dz_\alpha}{z_\alpha}\wedge dz_\beta\wedge d\ol z_J\otimes e_L
\ee
where now $\alpha, \beta$ in \eqref{mmp23} are multi-indices whose union is equal to $I$. Moreover $\alpha\subset \{1,\dots q\}$ and $\beta\subset \{q+1,\dots, n\}$. The length of $I$ and $J$ is $n-i$ and $i$, respectively.

The important remark is that the absolute value of the coefficients $\xi_{I\ol J}(\cdot, t)$ of the restriction \eqref{mmp23} are bounded uniformly with respect to $t\in V\setminus \Delta$. This is quickly seen thanks to the relation
\eqref{mmp10} which expresses the logarithmic form $\displaystyle \frac{dz_{n+m}}{z_{n+m}}$ as linear combination of $\displaystyle \frac{dz_{j}}{z_{j}}$ and the fact that the quotients $|z_{n+m}|/|z_j|$ are bounded from above on $U_{n+m}$. In particular,
we have
\be\label{rev2}
d\ol z_{n+m}= \sum_{j=1}^q\mu_jd\ol z_j
\ee
where the $\mu_j$ are uniformly bounded. 

\smallskip

\noindent In general, given two K\"ahler metrics $g_1$ and $g_2$ such that
$g_1\geq g_2$ then for any form $\gamma$ of type $(n-i, i)$ we have
\be\label{mp63}
|\gamma|_{g_1}\leq |\gamma|_{g_2} 
\ee
and in our setting this implies the inequality
\be\label{mmp24}
\left|\wt u_{|X_t}\right|_{\theta_\ep, h_{L,\ep}}^2\leq C\frac{e^{-\varphi_{L, \ep}}}{\ep^2}\sum_{I, J} |\xi_{I{\ol J}}|^2
\prod_{j\in \alpha} \log^2
\left|z_j\right|^2
\ee
at each point of the intersection $U_{n+m}\cap X_t$, cf. \eqref{mmp15}. This explains
the reason why we need to work with a metric having Poincar\'e singularities.

Next we remark that the local weights of the metric $h_L$ are bounded from below,
given the curvature hypothesis; this is also true for the perturbed metric $h_{L, \ep}$ (and actually the lower bound is independent of $\ep$). 

On the other hand, for each $\ep> 0$ there exists a constant $C_\ep> 0$ such that
we have 
\be\label{mmp25}
|t_m|^{2\ep}\prod_{j\in \alpha} \log^2
\left|z_j\right|^2\leq C_\ep
\ee
for any $z\in X_t\cap U$, cf. \eqref{mmp2}. We thus obtain 
\be\label{mmp26}
|t_m|^{2\ep}
\left|\wt u_{|X_t}\right|_{\theta_\ep, h_{L,\ep}}^2\leq C_\ep(\xi)
\ee
where $C_\ep(\xi)$ is uniform with respect to the point $t\in V\setminus \Delta$, as a consequence of 
\eqref{mmp24} and \eqref{mmp25}.

\noindent All in all, we infer the existence of a constant $C_\ep(\xi)$ such that we have
\be\label{mmp27}
|t_m|^{2\ep}
\int_{X_t}\left|\wt u_{|X_t}\right|_{\theta_\ep, h_{L,\ep}}^2\theta_{\ep}^n\leq C_\ep(\xi)
\ee
because the volume of each fiber $(X_t, \theta_{\ep})$ is bounded from above by a positive constant independent of $t$ (and $\ep$). Therefore the claim
is proved, and so is Theorem \ref{kernels, I} in case the line bundle $(L, h_L)$ verifies the hypothesis 
${\mathcal H}_2$. 
\medskip

\noindent If the curvature of $L$ is trivial when restricted to fibers of $p$, then we can use an arbitrary
K\"ahler metric on $X$ to define the norm on $\K^i$. Given this, the proof is completely identical, and we will not provide further details.
\end{proof}

\begin{remark}
  One can also prove the claim by a slightly different argument,
  which avoids the use of the restriction \eqref{mmp23}. The observation is that we have
\be\label{corrmp100}
\left|\wt u_{|X_t}\right|_{\theta_\ep, h_{L,\ep}}^2\leq \left|\wt u\right|_{\theta_\ep, h_{L,\ep}}^2
\ee
at each point of the fiber $X_t$. This follows by a simple linear algebra calculation. The right hand side of \eqref{corrmp100}
can be bounded from above by using the fact that $\theta_\ep$ is a metric with Poincar\'e singularities when restricted to $U$. In the inequality \eqref{mmp19} however, we have to work with the left hand side term of \eqref{corrmp100} directly, mainly because we have to obtain an upper bound which is independent of $\varepsilon$,
and the metric
$h_L$ is only assumed to be strictly positively curved in the fibers directions.
\end{remark}  

\begin{remark}
  The proof of Theorem \ref{kernels, I} is much easier than the one 
  establishing the \emph{positivity} of the direct image
  sheaf
  $$\F:= p_\star(K_{X/Y}+ L)$$
  in case $(L, h_L)$ has is semi-positively curved. The reason is that
  in order to show that the natural metric of $\F$ extends across the eventual singularities of the map $p$ 
  one has to show that the function
  $$\displaystyle t\to \Vert [u_t]\Vert_t^2$$
  is bounded \emph{from below} by a strictly positive constant, as soon as
  we have $\displaystyle [u_t]_{t_0}\not\in m_{0}\F_{t_0}$ where $t_0\in V\cap \Delta$ and $m_{0}\subset \O_{Y, t_0}$ is the maximal ideal associated to this point.
  
  \noindent On the other hand,
  in Theorem \ref{kernels, I} we assume that $L$ is endowed with a metric whose curvature satisfies very strong uniformity
requirements (with respect to the fibers of $p$).
As we will see in Theorem \ref{kernels},
  it is not always very easy to construct such metrics e.g. if $L=
  -K_{X/Y}$.
\end{remark}
\medskip

\noindent As already mentioned, the proof of Proposition \ref{ks} is
is practically contained in the arguments invoked for Theorem \ref{kernels, I}.
The only additional tool would be the following statement, generalizing
\cite{Del}.

\begin{lma}\label{mp6} Let $\Omega\subset Y$ be a small coordinate set
  centered at the point $y_0$. For each $i=1,\dots, m$ there exists
  a vector field $v_i$ defined on $p^{-1}(\Omega)$ with the
  following properties.
\begin{enumerate}

\item[(i)] We have $\displaystyle dp(v_i)= \frac{\partial}{\partial t_i}$ for $i=1,\dots, m-1$ and  $\displaystyle dp(v_m)= t_m\frac{\partial}{\partial t_m}$ on $\Omega$.
  \smallskip

\item[(ii)] On the open set $U\subset X$ as above we have

\be\label{mp10}
v_i= \frac{\partial}{\partial z_{n+i}}+
\sum_{j= 1}^q\psi^j_i\left(b_{m+n}z_j\frac{\partial}{\partial z_j}-
{b_jz_{m+n}}\frac{\partial}{\partial z_{m+n}}\right)+ 
\sum_{l= q+1}^n\psi^l_i\frac{\partial}{\partial z_l}
\ee
for $i=1,\dots, m-1$, well as
\be\label{mp11}
v_m= \frac{1}{b_{m+n}}z_{m+n}\frac{\partial}{\partial z_{n+m}}+ 
\sum_{j= 1}^q\psi^j_m\left(b_{m+n}z_j\frac{\partial}{\partial z_j}-
{b_jz_{m+n}}\frac{\partial}{\partial z_{m+n}}\right)+ 
\sum_{l= q+1}^n\psi^l_m\frac{\partial}{\partial z_l},
\ee
where the functions $\psi^i_j$
are smooth on $U$.
\end{enumerate}
\end{lma}

\begin{proof}
  The proof is straightforward: locally the lifting in (i) and (ii)
  clearly exist, given the structure of the map $p$. We glue them by a partition of unit; the logarithmic tangent bundle being invariant with respect to change of coordinates, the conclusion follows.   
\end{proof}
\smallskip

\noindent We remark that in order to construct the vector fields as in Lemma
\ref{mp6} it is enough to assume the following:
\begin{enumerate}
\item[{($s_1$)}] The singular values of the map $p$ are contained in a snc divisor $\Delta$;
  \smallskip
\item[{($s_2$)}] The inverse image $p^{-1}(\Delta)$ is a divisor
$W$ with snc support.    
\end{enumerate}

\medskip

\subsection{Proof of Theorem \ref{kernels, II}.} We will follow the method of Viehweg-Zuo, cf. \cite{VZ03}.

\noindent In order to simplify the writing, we introduce the following notation
\be\label{mp76}
E^s_{X/Y}:= {\mathcal R}^sp_\star\left(\Omega_{X/Y}^{n-s}\langle W\rangle\otimes
L\right)
\ee
for $s\geq 0$.

\noindent Then for each $s\geq 0$ we have the holomorphic map 
\be\label{mp69}
\tau^s: E^s_{X/Y}\to E^{s+1}_{X/Y}\otimes \Omega^1_Y\langle \Delta\rangle
\ee
as recalled at the beginning of the current section.
By hypothesis, there exists a section $\sigma$ of the bundle 
\be\label{mmp40}
 K_{X/Y}+ W- p^\star(\Delta)+ L
\ee
such that the map 
\be\label{mp71}
ks^{(1)}_\sigma: {\mathcal O}_Y\to E^1_{X/Y}\otimes \Omega^1_Y\langle \Delta\rangle.
\ee
given by the cup-product of $\sigma$ with the Kodaira-Spencer class is non-identically zero on
a non-empty open subset of $Y$.

\smallskip

\noindent Then the proof of Theorem \ref{kernels, II} is obtained as follows. Let
\be\label{mp81}
ks^{(j)}_\sigma: \O_Y\to
E^{j}_{X/Y}\otimes \Sym^{j}\Omega^1_Y\langle \Delta\rangle
\ee
be the map deduced from the iterated Kodaira-Spencer map
\eqref{mmp58}. 

We denote by $\Pi_j$ the projection 
\be\label{higher81}
E^{j}_{X/Y}\otimes \Sym^{j}\Omega^1_Y\langle \Delta\rangle\to \left(E^{j}_{X/Y}/\K^j\right)\otimes \Sym^{j}\Omega^1_Y\langle \Delta\rangle,
\ee
and let $i$ be the smallest integer such that
$\displaystyle \Pi_{i}\circ ks^{(i)}_\sigma$
is identically zero. This means that the image of $ks^{(i)}_\sigma$ is contained in
\be\label{mp83}\displaystyle {\mathcal K}^{i}\otimes\Sym^{i}\Omega^1_Y\langle \Delta\rangle
\ee
and that $ks^{(i)}_\sigma$ is not identically zero.
Indeed, if $ks^{(i)}_\sigma$ would be zero, then the image of $ks^{(i-1)}_\sigma$ is contained in $\displaystyle {\mathcal K}^{i-1}\otimes\Sym^{i-1}\Omega^1_Y\langle \Delta\rangle$. This contradicts our choice of $i$.

Note that we have $i\geq 1$ thanks to the assumption $({\mathcal H})_3$. We obtain a non-identically zero section of the bundle \eqref{mp83}. This in turn defines a non-trivial map
\be\label{mp84} {\mathcal K}_f^{i \star}\to \Sym^{i}\Omega^1_Y\langle \Delta\rangle.
\ee
since the dual of $\K^i$ is the same as the dual of $\K^i_f$.
We remark that this already shows that the bundle $\Sym^{i}\Omega^1_Y\langle \Delta\rangle$
has a semi-positively curved subsheaf (i.e. the image of the map \eqref{mp84}).
\smallskip

\noindent If the curvature of $(L, h_L)$ satisfies the hypothesis
$({\mathcal H}_2)$ and if moreover there exists an open subset $\Omega\subset Y$ together with a positive $\varepsilon_0> 0$ such that
$$\Theta_{h_L}(L)\leq -\varepsilon_0p^\star (\omega_Y)$$
on $p^{-1}(\Omega)$, then the curvature properties of the kernels
$\K^is$ are considerably better:
\smallskip

\noindent $\bullet$ The sheaf $\K^i$ is admits a singular semi-negatively curved
Hermitian metric which is smooth on a Zariski open subset of $Y$.
\smallskip

\noindent $\bullet$ For any local holomorphic section $[u_t]$ of the bundle $\displaystyle \K^i|_{\Omega}$ we have
\be\label{mmp50} 
dd^c\log \Vert [u]\Vert^2\geq \ep_0\sum_{j=1}^m \sqrt{-1}dt_j\wedge d\ol t_j. 
\ee
\medskip

\noindent In this context we recall the following result, which concludes the proof of the last part of Theorem \ref{kernels, II}. 

Let $Y$ be a smooth projective variety, and let $\F$ be a torsion free coherent sheaf on $Y$. We consider 
${\mathbb P}(\F):= \Proj (\bigoplus_{m \ge 0} S^m(\F))$ the scheme over $Y$ associated to $\F$, together with the projection $\pi : {\mathbb P}(\F) \to Y$.
We denote by $\O_{\F}(1)$ the tautological line bundle on ${\mathbb P}(\F)$. Let $Y_1 \subset Y$ be a Zariski open subset on which $\F$ is locally free (in particular ${\mathbb P}(\F)$ is smooth over $Y_1$) and
${\rm codim}_{Y}(Y\setminus Y_{1}) \ge 2$.
\medskip

\noindent The following statement is established in \cite{PT}.

\begin{thm} \label{bigness}\cite{PT}
We suppose that $\O_{\F}(1)|_{\pi^{-1}(Y_{1})}$ admits a singular Hermitian metric $g$ with semi-positive curvature, and 
that there exists a point $y \in Y_1$ such that
${\mathcal I}(g^k|_{{\mathbb P}(\F_y)})=\O_{{\mathbb P}(\F_y)}$ for any $k>0$, where ${\mathbb P}(\F_{y})=\pi^{-1}(y)$. Assume moreover that there exists an open neighborhood $\Omega$ of $y$ such that $\Theta_{g}\left(\O(1)\right) - \pi^{\star}\omega_Y \ge 0$ on $\pi^{-1}(\Omega)$.
Then $\F$ is big.
\end{thm}

\medskip

\noindent This result together with the properties of the kernels
$\K^j_f$ summarized in the bullets above show that each of $\K^{j\star}_f$ is
big; in particular this applies to $\K^{i\star}_f$. Let $A$ be an ample line bundle.
We obtain a section of $\Sym^M\K^{i\star}_f\otimes A^{-1}$ for $M\gg 0$, which induces a section of $\otimes^{iM}\Omega^1_Y\langle \Delta\rangle\otimes A^{-1}$ thanks to \eqref{mp84}.
\qed

\begin{remark}
  The hypothesis $({\mathcal H}_i)$ (for $i=1,2$) together with the existence of
  a section $\sigma$ of $K_{X/Y}+ W- p^\star(\Delta)+ L$ have a few consequences on the map $p$ itself which we will now discuss. We write
\be\label{higher100}
K_{X/Y}= L_1- L
\ee  
where $\displaystyle L_1:= \big(K_{X/Y}+ W- p^\star(\Delta)+ L\big)+ p^\star(\Delta)- W$. In particular, $L_1$ is an effective line bundle. As for
$-L$, it is semi-positive and strictly positive/trivial on the fibers of
$p$. Therefore, the relative canonical bundle $K_{X/Y}$ is relatively big in the first case, and pseudo-effective in the second.

If moreover the curvature of $L$ verifies the requirements in
the second part of Theorem \ref{kernels, II}, then for any $m$ large enough the bundle $\det p_\star(mK_{X/Y})$
is big. From this perspective, Theorem \ref{kernels, II} is similar to (but slightly weaker than) the main statement of Popa-Schnell in \cite{PoSch}. 
\end{remark}

\subsection{Metrics on the relative canonical bundle}

\noindent In this section we will establish the following result; afterwards we will show that it implies Theorem \ref{kernels}.

\begin{thm}\label{mp5}
  Let $p:X\to Y$ be an algebraic fiber space, and let $\omega$ be a fixed reference metric on $X$. We assume that the following properties are
  satisfied.
\begin{enumerate}

\item[(1)] The canonical bundle of the generic fiber of $p$ is ample.
\smallskip  

\item[(2)] There exists a positive $m\gg 0$ such that the line bundle
$$\det p_\star(mK_{X/Y})$$   
is big.  
\end{enumerate}    

\noindent Then there exist an effective, $p$-contractible divisor $\Xi$ on $X$ and a singular
metric $\displaystyle h_{X/Y}= e^{-\varphi_{X/Y}}$ on $K_{X/Y}+ \Xi$ such
that:

\begin{enumerate}

\item[(i)] There exists $\varepsilon_0> 0$ for which we have $dd^c\varphi_{X/Y}\geq \varepsilon_0\omega$ in the sense of currents on $X$, i.e. the curvature of $(K_{X/Y}+ \Xi, h_{X/Y})$ is strongly positive.
\smallskip

\item[(ii)] The restriction $\displaystyle h_{X/Y}|_{X_y}$ is non-singular for each
  $y\in B$, where $B$ is a Zariski open subset of $Y$ (which can be described in a precise manner).
\end{enumerate}
\end{thm}

\noindent Before presenting the arguments of the proof, a few comments about the preceding statement.


\begin{remark}
  If instead of the hypothesis (1) above we assume that $\displaystyle K_{X_y}$ is big, then the metric $e^{-\varphi_{X/Y}}$ can be constructed so that it has the same singularities as the \emph{metric with minimal singularities} on this bundle. It would be really interesting to develop the techniques in sections (1)-(6) under the
  assumption that the bundle $L$ is endowed with a sequence of smooth metrics $h_\ep$ such that the curvature $\theta_\ep$ admits a lower bound
  $$\theta_\ep \geq (\varepsilon_0- \lambda_\ep)\omega$$
where $\varepsilon_0> 0$ and $\lambda_\ep$ is a family of positive functions uniformly bounded from above and converging to zero almost everywhere. 
\end{remark}

\begin{remark}
  If the Kodaira dimension of the fibers of $p$ is not maximal we can
  derive a similar result --of course,
  the conclusion (i) has to be modified accordingly.
  However, in the absence of the hypothesis (2) it is not clear what we can expect.
\end{remark}
  
\begin{proof}
  The arguments which will follow rely on two results. We recall next the
  first one, cf. \cite{CP}.


\begin{thm}\label{mp31}
  \cite{CP} We assume that $K_X$ is ${\mathbb Q}$-effective when restricted to the generic fiber of $p$. 
  For any positive integer $m$ there exists a real number $\eta_0> 0$ and a divisor $\Xi$ on $X$ such that the codimension of $p(\Xi)$ is at least two, and such that
  the difference
\begin{equation}\label{mp2}
K_{X/Y}+ \Xi - \eta_0p^\star\left(\det\big(p_\star(mK_{X/Y})\big)\right)
\end{equation}  
is pseudo-effective.  
\end{thm}  

\noindent Actually what we will use is not this result in itself, but its proof
--this will be needed in order to establish (ii) above.
Hence the plan for the remaining part of this subsection is to review the main
parts of the proof as in \cite{CP} (see also the references therein), and to extract the statement we want --this concerns the singularities of the current in \eqref{mp2}. 
\medskip

\noindent $\bullet$ If the map $p$ is a smooth submersion, then the arguments in \cite{CP} are borrowed from Viehweg's weak semi-stability theorem proof, cf. \cite{Vbook}. The first observation is that given any vector bundle 
$E$ of rank $r$ we have a canonical injection
\be\label{mp32}
\det(E)\to \otimes^rE,
\ee
and hence a nowhere vanishing section of the bundle
\be\label{mp33}
\otimes^rE\otimes \left(\det(E)\right)^\star.
\ee
We apply this to the direct image
\be\label{mp34}
E:= p_\star\left(m_0K_{X/Y}\right)
\ee
which is indeed a vector bundle, by the invariance of plurigenera \cite{Siupg}. Here $m$ is a positive integer, large enough so that the multiple $\displaystyle m_0K_{X_t}$ is very ample, and $r$ is the rank of the
direct image \eqref{mp34}.

In this case the bundle $\otimes^rE$
can be interpreted as direct image of the relative pluricanonical bundle of
the map   
\be\label{mp35}
p^{(r)}:X^{(r)}\to Y
\ee
where $\displaystyle X^{(r)}:= X\times_Y\dots\times_Y X$ is the $r^{\rm th}$ fibered product corresponding to the map $p:X\to Y$. Note that $X^{(r)}\subset \times^rX$ is a non-singular submanifold of the $r$-fold product of $X$ with itself (thanks to the assumption that $p$ is
a submersion).

The important observation is that we have the formulas
\be\label{mp36}
p^{(r)}_\star\left(m_0K_{X^{(r)}/Y}\right)= \otimes^rp_\star\left(m_0K_{X/Y}\right)
\ee
as well as
\be\label{mp37}
K_{X^{(r)}/Y}= \prod \pi_i^\star(K_{X/Y})
\ee
where $\pi_i:X^{(r)}\to X$ is induced by the projection on the $i^{\rm th}$ factor. 

By relations \eqref{mp33} and \eqref{mp36} we have a section, say
$\sigma$, of the bundle
\be\label{mp38}
m_0K_{X^{(r)}/Y}-p^{(r)\star}\left(\det p_\star\left(m_0K_{X/Y}\right) \right).
\ee
The formula \eqref{mp37} shows that the restriction of the bundle $K_{X^{(r)}/Y}$ to the diagonal $X\subset X^{(r)}$ is precisely $rK_{X/Y}$, so all in all the restriction of the bundle in \eqref{mp38} to $X$ coincides with \eqref{mp2}. However, we cannot use directly the section $\sigma$ in order to conclude, since by construction we have $\displaystyle \sigma|_{X}\equiv 0$.
\smallskip

To bypass this difficulty, we will use the following superb trick invented by E. Viehweg, cf. \cite{Vbook}. Let $\ep_0> 0$ be a small enough positive rational number such that the pair
\be\label{mp39}
\left(X^{(r)}, B\right)
\ee
is klt, where the boundary $B:= \ep_0Z_\sigma$ is the $\ep_0$ multiple of the
divisor $Z_\sigma:= (\sigma= 0)$.

The next step is to apply the results in \cite{BP} and deduce that there exists a fixed ample line bundle $A_Y$ on $Y$ such that for any
$k\geq 1$ divisible enough
the restriction map
\be\label{mp40}
H^0\left(X^{(r)}, k\left(K_{X^{(r)}/Y}+ B\right)+ p^{(r)\star}A_Y\right)
\to H^0\big(X^{(r)}_y, k(K_{X^{(r)}_y}+ B|_{X^{(r)}_y})\big)
\ee
is surjective, for any $y\in Y\setminus \Delta$. Now the bundle
on the right hand side of \eqref{mp40} is equal to
\be
k(1+ \ep_0m_0)\prod_{i=1}^r \pi_i^\star \left(K_{X_y}\right), 
\ee
so it has many sections, given that $K_{X_y}$ is ample.
Moreover, all the sections in question are automatically
$L^{2/k}$ integrable, given that $(X^{(r)}, B)$ is a klt pair, and so it is its restriction to a generic fiber $X^{(r)}_y$ of our
map $p^{(r)}$. For the complete argument concerning the surjectivity of the map
\eqref{mp40} we are referring to \cite{CP}
(and the references therein).

In particular, for each $k$ divisible enough the bundle
\be\label{mp41}
k\left(K_{X^{(r)}/Y}+ B\right)+ p^{(r)\star}A_Y
\ee
admits a holomorphic section whose restriction to the diagonal
$X\subset X^{(r)}$ is not identically zero. Indeed, we have (many)
sections of the bundle \eqref{mp40} on
fiber $X_y\times\dots\times X_y$ whose restriction to the diagonal
$X_y\subset X_y\times\dots\times X_y$ is not identically zero.

\smallskip

\noindent In conclusion, this shows the existence of a positively-curved metric on the bundle
\be\label{mp42}
r(1+\ep_0m_0)K_{X/Y}-\ep_0p^\star\det \big(p_\star(m_0K_{X/Y})\big)+ \frac{1}{k}p^\star A_Y
\ee
whose restriction to the generic fiber $X_y$ of $p$ is induced by the
sections of multiples of $\displaystyle K_{X_y}$.
\medskip

\noindent $\bullet$ The general case where $p$ is simply a surjective map is much more involved; the complications are arising from the fact that the fibered product $X^{(r)}$ is no longer smooth. However, as shown in \cite{CP} one obtains a similar result, modulo adding the divisor $\Xi$ which projects in codimension greater than two.
We will not reproduce here the proof, but just mention
that modulo the singularities of $X^{(r)}$, the structure of the argument is identical to the case detailed at the preceding bullet --in
particular, we control the singularities of the restriction of the
metric on \eqref{mp42} to the generic fiber $X_y$.

\begin{remark}
  It would be really interesting to have a more direct 
  proof of Theorem \ref{mp31}, i.e. without using Viehweg's trick.
\end{remark}  
\bigskip

\noindent Recall that
we have an integer $m_0$ above such that the following conditions are satisfied.
\begin{enumerate}

\item[(i)] The bundle $\displaystyle m_0K_{X_y}$ is very ample.
\smallskip

\item[(ii)] The bundle $\displaystyle \det \left(p_\star(m_0K_{X/Y})\right)$ is big.
\end{enumerate}  

\noindent By the point (ii) above, there exists a positive
integer $m_1\gg 0$ such that
\be\label{mp47}
m_1\det \left(p_\star(m_0K_{X/Y})\right)\simeq A_Y+ E_Y
\ee
where $E_Y$ is an effective divisor on $Y$.

The properties of the bundle \eqref{mp42} combined with our previous considerations show the existence of a metric $e^{-\psi_{X/Y}}$ on the bundle
$K_{X/Y}+ \Xi$ with the following properties.
\begin{enumerate}

\item[(a)] The metric $e^{-\psi_{X/Y}}$ is semi-positively curved, and it has algebraic singularities (meaning that its local weights are of the form log of a sum of squares of holomorphic functions, modulo a smooth function). Moreover, the restriction
$\displaystyle e^{-\psi_{X/Y}}|_{X_y}$ is smooth, for any $y$ belonging to a Zariski open subset of $Y$.
  \smallskip

\item[(b)] The form $dd^c\psi_{X/Y}$ is smooth and definite positive on the $p$-inverse image of a (non-empty) Zariski open subset. 

\end{enumerate}
\medskip

\noindent We see that the main difference between the properties
(a) and (b) of the metric $e^{-\psi_{X/Y}}$ and the conclusion of Theorem \ref{mp5} is strict positivity, i.e. the existence of $\varepsilon_0$ as in (i). In order to finish the proof, we will use a result due to
M. Nakamaye, cf. \cite{Nak} concerning the augmented base loci of nef and big line bundles.

As a preparation for this, we blow-up the ideal corresponding to the
singularities of $e^{-\psi_{X/Y}}$; let $\pi: \wt X\to X$ be the
associated map. We write
\be\label{frei1}
\pi^{\star}\left(K_{X/Y}+ \Xi\right)\simeq {\mathcal L_1}+ {\mathcal L_2}
\ee
such that ${\mathcal L_j}$ above are ${\mathbb Q}$-line bundles induced by the decomposition of the inverse image of the curvature of $e^{-\psi_{X/Y}}$
into smooth and singular parts denoted by $\theta_1$ and $\theta_2$,
respectively. By the properties (a) and (b) above, $\theta_1$ is smooth and definite positive in the complement of an algebraic set which projects into a strict subset of $Y$. As for $\theta_2$, it is simply the current of integration on an algebraic subset of $\widetilde X$ which equally projects properly into $Y$.
\smallskip

\noindent In particular, the bundle ${\mathcal L_1}$ is nef and big (given that it admits a metric whose curvature is $\theta_1$).
We denote by $\displaystyle {\rm B}_+({\mathcal L_1})$ the stable base locus of the ${\mathbb Q}$-bundle ${\mathcal L_1}- A$, where
$A$ is any small enough ample on $\widetilde X$.
We recall the following result, cf. \cite{Nak}.

\begin{thm}\cite{Nak} The algebraic set  
$\displaystyle {\rm B}_+({\mathcal L_1})$ is the union of $W\subset {\widetilde X}$ such that $\displaystyle \int_W\theta_1^d= 0$,
where $d$ is the dimension of $W$.  
\end{thm}  
\noindent Given the positivity properties of the form $\theta_1$, we infer that the
set $p\circ \pi \left({\rm B}_+({\mathcal L_1})\right)$ is strictly contained in $Y$. 
This then implies that we can endow ${\mathcal L_1}$ with a metric whose curvature can be written as the sum of a K\"ahler metric (given by the ample $A$ above) and a closed positive current whose singularities are $p\circ \pi$-vertical. By combining it with the metric on ${\mathcal L_2}$ we get a metric with the same properties on the
inverse image $\pi^\star(K_{X/Y}+ \Xi)$. 

\smallskip

\noindent The metric $h_{X/Y}$ is obtained by push-forward, and
Theorem \ref{mp5} is proved.
\end{proof}

\medskip

\subsection{Proof of Theorem \ref{kernels}} The hypothesis of \ref{kernels} together with the results in \cite{Kol87}, \cite{Kawa}, \cite{PoSch} that for each
$m\gg 0$ the bundle
\be\label{mmp60}
\det\left(p_\star(mK_{X/Y})\right)
\ee
is big.
This can equally be obtained along the line of arguments in this paper as follows.
It is shown in \cite{Berndtsson2} that the curvature form of the vector 
bundle $p_\star(mK_{X/Y})$ has no zero eigenvalue --as if not, the ``maximal variation" hypothesis would be contradicted. Thus the curvature form of its determinant is semi-positive, and strictly positive at some point.
The fact that it is big follows e.g. by holomorphic Morse inequalities.
We now apply Theorem \ref{mp5} and \ref{kernels} is proved.

\subsection{Proof of Theorem \ref{VZ}} This is almost a linear combination of the results obtained in the previous sections: we define the
$L$ as follows

\be\label{031005mp}
L:= -K_{X/Y}+ \sum(t^i-1)W_i
\ee
so that the bundle \eqref{mmp40} is trivial (in particular, it has a section).
\smallskip

\noindent We intend to use Theorem \ref{kernels, II} in order to conclude, but 
we do not seem to prove that the bundle
\eqref{031005mp} admits a metric satisfying the hypothesis
$\left({\mathcal H}_2\right)$. However, in Corollary \ref{approx} we construct
a family of metrics whose curvature properties represent an
approximation of this hypothesis. We will show next that we can adapt the proof of Theorem \ref{kernels, II}, and obtain the same conclusion by using Corollary \ref{approx} instead of $\left({\mathcal H}_2 \right)$.

We will denote by $h_L$ the metric on \eqref{031005mp} induced by the metric
$h_{X/Y}^{(1)}$. Note that the curvature current
corresponding to this metric is semi-positive on $X$, and strictly positive
on generic fibers of $p$ by Theorem \ref{bo_powers}.

Let $[u_t]$ be a local holomorphic section of the bundle $\K^i$, defined on 
the open co-ordinate subset $V\subset Y$ centered at a generic point of $\Delta$.
We have to show that the quantity
\be\label{031101mp}
\sup_{t\in V\setminus \Delta}\log\Vert [u_t]\Vert^2_t< \infty
\ee
is finite. To this end, we follow the same path as in the proof of \ref{kernels, I}, so we will only highlight the main differences next.

By Corollary
\ref{approx} we have 
\be\label{031102mp}
h_L|_{X_t}= \lim_{\eta\to 0}h_{X/Y}^{(\eta)}|_{X_t}
\ee
for each $t\in V\setminus \Delta$. For each $0<\ep\ll \eta$, we define the
perturbation $\displaystyle h_{L, (\ep,\eta)}$ of the metric $h_{X/Y}^{(\eta)}$ precisely as in
\eqref{mmp16} (meaning that $h_L$ in \eqref{mmp16} is replaced with $h_{X/Y}^{(\eta)}$ in the actual context). Then we have
\be\label{031103mp}
\sup_{t\in V\setminus \Delta}|t_m|^{C(1-\eta)}\Vert [u_t]\Vert^2_{t, (\ep, \eta)}\leq
C(\ep, u)< \infty,
\ee
where $C$ in \eqref{031103mp} is a fixed constant, large enough compared with the
multiplicities $b^i$ in the expression of the divisor $p^{-1}(\Delta)$. The inequality  \eqref{031103mp} above is
established by the same procedure as the proof of the
Claim --the only slight difference here is the presence of a factor of
order $\O(1-\eta)$ with the wrong sign in the metric $h_{X/Y}^{(\eta)\star}$, accounting for 
the first negative term in the curvature estimate (a) in Corollary \ref{approx}.
This term is tamed by the factor $|t_m|^{C(1-\eta)}$.
The inequality \eqref{031103mp} together with the mean inequality
show that \eqref{031101mp} hold true.
\smallskip

\noindent Hence, under the hypothesis of Theorem \ref{VZ} the conclusion of
theorems \ref{kernels, I} and \ref{kernels, II}, respectively hold true
provided that the bundle $(L, h_L)$ is chosen as in \eqref{031005mp} above.

\subsection{Proof of Theorem \ref{CY17}} Let $p:X\to Y$ be a Calabi-Yau family, which has maximal variation. By the results in \cite{Kol87}, \cite{Kawa}, \cite{PoSch} we infer that for all
$m\gg 0$ divisible enough the bundle
\be\label{mmp100}
\L:= \left(p_\star(mK_{X/Y})\right)^{\star\star}
\ee
is big. Actually the metric version of this statement is true, as we will see after recalling a few facts.

Let $y\in Y$ be a regular value of the map $p$. Since $c_1(X_y)= 0$, there
exists a positive integer $m$ such that the bundle $\displaystyle mK_{X_y}$
admits a nowhere vanishing section say $s_y$. The invariance of plurigenera \cite{Siupg} shows that for any co-ordinate open subset $V$ containing $y$
there exists a section $s$ of the bundle $\displaystyle mK_{X/Y}|_{p^{-1}(V)}$
whose restriction to the fiber $X_y$ equals $s_y$.

By shrinking $V$, we can assume that $s$ is nowhere vanishing, and therefore it gives a trivialization of the bundle $\L|_V$. With respect to this trivialization, the local weight of the metric on $\L$ in \cite{Berndtsson1},
\cite{PT} is given by the function
\be\label{mmp101}
\varphi_\L(w)= \log\int_{X_w}|s|^{2/m}
\ee
where the expression under the integral in \eqref{mmp101} is the volume element on $X_w$ induced by the restriction of $s$.    
\medskip

\noindent It is proved in \cite{PT} that the metric \eqref{mmp101} is semi-positively curved. Under the hypothesis of \ref{CY17} the following stronger result hold true.

\begin{lma}\label{cymetric}
We assume that the Calabi-Yau family $p$ has maximal variation.
Then the curvature of $(\L, e^{-\varphi_\L})$ is positive definite on each compact subset of a non-empty Zariski open subset of $Y$. In particular, $\L$ is big.   
\end{lma}  
\begin{proof} We observe that the metric $e^{\varphi_\L}$ is smooth on a Zariski open subset $Y_0\subset Y$. Let $\Theta_\L$ be the corresponding curvature
  form. We claim that $\Theta_\L^{\dim (Y)}> 0$ on $Y_0\cap U$, where $U\subset Y$ is such that the restriction of the  Kodaira-Spencer map $ks|_U$ to $U$
  is injective. The argument goes as follows: if we have $\Theta_\L^{\dim (Y)}= 0$, then the kernel of $\Theta_\L$ would define a foliation which is not necessarily holomorphic, but whose leaves are holomorphic. The fact that the leaves are holomorphic is basically a consequence of the fact that $\Theta_\L$ is of (1,1)--type.
  We refer to \cite{Kol87} and the references therein for a complete explanation.
  
  We consider a holomorphic disk $\D$ contained in a generic leaf of this (local) foliation. Then the family $f: X_\D\to \D$ induced by $p$ a submersion, where $X_\D:= p^{-1}(\D)$, and moreover the curvature of the direct image $f_\star(mK_{X_\D/\D})$ is equal to zero. If $m=1$ we invoke the results in \cite{Choi} 
to infer that this forces the
vanishing of the Kodaira-Spencer class of $f$. If $m\geq 2$, we argue as follows. By the definition of the
metric $\varphi_{X/Y}$ in Section 8.2 together with the fact that the curvature of $f_\star(mK_{X_\D/\D})$
is zero, we infer that
$$i\ddbar \varphi_{X/Y}= 0.$$ 
This implies the vanishing of $\omega_{WP}$, cf. \eqref{WP100}. 
We therefore obtain a contradiction, so Lemma \ref{cymetric} is proved.
\end{proof}
\smallskip

\noindent As a consequence, we obtain the following statement.

\begin{cor}\label{cymetricII}
  Let $p:X\to Y$ be a Calabi-Yau family which has maximal variation. Then the
  relative canonical bundle $K_{X/Y}$ has the following property.
 There exists a positive integer $m$ such that the curvature current $\Theta$ of the $m$-Bergman metric
  $\displaystyle e^{-\varphi^{(m)}_{X/Y}}$ is semi-positive on $X$ and greater than the inverse image
  of a metric on the pre-image of a non-empty open subset of $Y$. Moreover, we have
\be\label{mmp102}
\Theta\geq \sum_i(t^i- 1)[W_i]  
\ee
on $X$. 
\end{cor}

\begin{proof}
  Indeed, in our set-up the $m$--Bergman metric has the following expression
\be\label{mmp103}
e^{\varphi^{(m)}_{X/Y}(x)}= \frac{|s|_x^{2/m}}{\int_{X_y}|s|^{2/m}}
\ee
where $y= p(x)$. Therefore we have $\Theta= p^\star(\Theta_\L)$ and
our statement follows from Lemma \ref{cymetric}.
\end{proof}

\noindent Theorem \ref{CY17} is a direct consequence of Theorem \ref{kernels, I}
together with Theorem \ref{kernels, II}.
This is the end.


\end{document}